\documentclass[notitlepage,a4,10.5pt]{article}

\usepackage[colorlinks=true,citecolor=blue,linkcolor=blue]{hyperref}

\usepackage[left=2.5cm,right=2.5cm,top=2.9cm,bottom=2.9cm]{geometry}

\usepackage{comment}
\usepackage{xurl}

\usepackage{tikz}
\usetikzlibrary{arrows.meta,positioning,calc,shapes,fit}
\usepackage{subcaption}

\usepackage{booktabs}

\usepackage{tabularx}
\usepackage{array}
\newcolumntype{Y}{>{\raggedright\arraybackslash}X}
\usepackage{booktabs}
\usepackage{multirow}

\usepackage{graphicx}

\usepackage{amsmath}
\numberwithin{equation}{section}

\usepackage{amsmath,array,booktabs,longtable}

\usepackage{amsfonts}
\usepackage{amssymb}

\usepackage{mathrsfs}
\usepackage{amsthm}
\usepackage{float}

\usepackage{authblk}

\usepackage{enumitem}

\usepackage{stmaryrd}

\usepackage{xcolor}
\usepackage{appendix}

\newtheorem{thm}{{Theorem}}[section]

\theoremstyle{definition}
\newtheorem{mydef}[thm]{Definition}

\theoremstyle{remark}
\newtheorem{remark}[thm]{\textbf{Remark}}

\newtheorem{corollary}[thm]{\textit{Corollary}}
\newtheorem{prop}[thm]{\textbf{Proposition}}

\newtheorem{lemma}[thm]{\textbf{Lemma}}

\theoremstyle{plain}
\newtheorem{conjecture}[thm]{Conjecture}

\definecolor{myyellow}{RGB}{255,255,180}

\begin{document}

\title{
An Analytic Toroidal MHD Equilibrium with Trivial Euclidean Symmetry Group
}

\author{
Kooper de Lacy\thanks{School of Mathematics and Statistics, University of Western Australia, Perth, WA, 6009, Australia   
\ Email: \href{kooper.delacy@research.uwa.edu.au}
{kooper.delacy@research.uwa.edu.au} 
}
\qquad 
Naoki Sato\thanks{National Institute for Fusion Science,  322-6 Oroshi-cho Toki-city, Gifu 509-5292, Japan,\ Email: \href{sato.naoki@nifs.ac.jp}{sato.naoki@nifs.ac.jp}
}  
{}\thanks{Graduate School of Frontier Sciences, The University of Tokyo,  
Kashiwa, Chiba 277-8561, Japan}
}

\date{\today}
\setcounter{Maxaffil}{0}
\renewcommand\Affilfont{\itshape\small}

    \maketitle

\begin{abstract}
    We provide a method to construct a real-analytic ideal MHD equilibrium $(B,p)$ in a solid toroidal domain $\Omega\subset\mathbb{R}^{3}$, with $B$ tangential to the boundary $\partial\Omega$. The regular pressure level sets foliate $\Omega\setminus\Gamma$ by nested tori surrounding a magnetic axis $\Gamma$, with $\nabla p\neq0$ away from the axis. We construct an example where the magnetic field is nonvanishing on $\Omega$, and the rotational transform equals 2 on every regular pressure surface. The equilibrium has a trivial Euclidean symmetry group. It is non-isolated in the sense that it belongs to a real-analytic one-parameter family connecting it to an axisymmetric equilibrium. The toroidal domain, magnetic field and pressure vary smoothly along this family. Under the MHD--Euler correspondence, the construction also yields a real-analytic steady incompressible Euler flow with a nonconstant Bernoulli function. This constructively disproves one formulation of Grad’s conjecture by exhibiting a smooth toroidal MHD equilibrium with nested pressure surfaces whose Euclidean symmetry group is trivial.
\end{abstract}


\section{Introduction}

The magnetohydrodynamic (MHD) equilibrium equations describe the equilibria of an ideal, perfectly conducting plasma.
Normalising the magnetic permeability to one, the governing equations are

\begin{equation}
    \nabla\cdot B=0 \, ,\qquad (\nabla\times B)\times B=\nabla p \, . \label{eq:mhs}
\end{equation} 
Here $B$ is a magnetic field and $p$ is a pressure on a solid torus $\Omega\subset\mathbb R^3$, by which we mean a compact domain diffeomorphic to $\mathbb D\times \mathbb S^1$ (where $\mathbb D$ is the closed unit disk, and $\mathbb S^1=\mathbb R/(2\pi\mathbb Z)$). We seek analytic functions $(B,p)$
on $\Omega$ that satisfy \eqref{eq:mhs}.

We consider equilibria whose pressure is constant on the boundary and the pressure level sets are nested tori that collapse onto a closed curve, the magnetic axis. Additionally, we restrict to the case where \(\nabla p \neq 0\) away from the magnetic axis. Taking the inner product of the second equation in \eqref{eq:mhs} with $B$ gives
\begin{equation}
    B\cdot \nabla p=0 \, . \label{eq:bc}
\end{equation}
Hence $B$ is tangent to each pressure level set, including the boundary.

Grad conjectured that smooth equilibria with nonconstant pressure exist only for symmetric configurations~\cite{Grad1967}. He later reformulated this conjecture as the nonexistence, apart from the stated symmetric exceptions, of families of toroidal equilibria with nested pressure surfaces depending smoothly on a parameter~\cite[p.~35]{Grad1985}. Constantin, Drivas, and Ginsberg formalised these statements as follows.

\begin{conjecture}[Grad's conjecture~\cite{Grad1967,Grad1985}, quoted up to notational differences from Constantin, Drivas, and Ginsberg~{\cite[Conjecture~1]{ConstantinDrivasGinsberg2021}}]\label{conj:grad}

    Any non-isolated and non-vanishing (away from the “magnetic axis”) smooth unforced MHD equilibrium on a (topologically toroidal or cylindrical) domain \(\Omega \subset \mathbb R^3\) that has a pressure possessing nested level sets which foliate \(\Omega\) has either plane-reflection, axial or helical symmetry.
\end{conjecture}

An equilibrium is \emph{isolated} if, in some suitable topology, every nearby equilibrium is a trivial rescaling or translation of it~\cite{ConstantinDrivasGinsberg2021}.


\subsection{Statement of the main results}

A \emph{symmetry} of $(\Omega,B)$ is a Euclidean isometry \(G\) such that $G(\Omega)=\Omega$ and
\begin{equation}
    G(x)=Ox+\tau \, , \qquad \qquad
    B(Gx)=\eta\,OB(x) \, ,\qquad \qquad (\text{for all } x\in\Omega) \, ,\label{eq:sym}
\end{equation}
with $O\in\mathrm O(3)$, $\tau\in\mathbb R^3$, and $\eta\in\{1,-1\}$. Reflections are included. The sign $\eta$ admits symmetries that reverse the field, such as stellarator symmetry~\cite{DewarHudson1998}, and accommodates the convention that $B$ is an axial vector. Force balance implies that every symmetry also preserves $p$ (Subsection~\ref{subsec:asymmetry}).

\vspace{.5em}

\begin{thm}[Existence]\label{thm:main}
For every $\varepsilon\in[0,1]$, there exists a solid torus $\Omega_\varepsilon$
with analytic boundary and a solution $(B,p)$ of
\eqref{eq:mhs}, real analytic on  $\Omega_\varepsilon$,
with the common magnetic axis
\[
 \Gamma=\big\{(\cos\theta,\sin\theta,0):\theta\in\mathbb S^1 \big\} \, ,
\]
and the following properties.
\begin{enumerate}
\item $B$ is nowhere zero and tangent to $\partial\Omega_\varepsilon$. \label{item:part1}
\item The critical set of $p$ in $\Omega_\varepsilon$ is exactly
$\Gamma$. Every other pressure level in $\Omega_\varepsilon$ is a torus,
and these tori form a nested family. \label{item:part2}
\item Every field line is closed, and every field line other than $\Gamma$ has rotational transform $2$. \label{item:part3}
\item If $\varepsilon>0$ and an isometry $G$ satisfies \eqref{eq:sym} for
$(\Omega_\varepsilon,B)$, then $G=\mathrm{Id}$ and $\eta=1$. \label{item:part4}
\item At $\varepsilon=0$, the domain and equilibrium are invariant under
all rotations about the vertical axis. \label{item:part5}
\end{enumerate}
\end{thm}

\vspace{.5em}

\begin{thm}[Regularity and non-isolation]\label{thm:nonisolation}
The family of Theorem~\ref{thm:main} can be chosen with embeddings
$X_{\mathbb D}:\mathbb D\times\mathbb S^1\longrightarrow\mathbb R^3$ whose images are
$\Omega_\varepsilon$, for $\varepsilon\in[0,1]$. The embedding $X_{\mathbb D}$ and
the compositions $B\circ X_{\mathbb D}$ and $p\circ X_{\mathbb D}$ are jointly real analytic
on 
$\mathbb D\times\mathbb S^1\times[0,1]$. Furthermore, every member of the family is non-isolated in the $C^\infty$ topology on the triples $(X_{\mathbb D},B\circ X_{\mathbb D},p\circ X_{\mathbb D})$. 
\end{thm}

\vspace{.5em}

Throughout the family, we write $X_{\mathbb D}$, $B$ and $p$ with their dependence on the parameter $\varepsilon$ suppressed.

Each member with $\varepsilon>0$ is therefore an MHD equilibrium without symmetry. It is joined to the axisymmetric member by a path along which the domain, the field and the pressure vary analytically, and every member on this path other than the endpoint has no symmetry. Such a member satisfies the hypotheses of Conjecture~\ref{conj:grad} and has none of the symmetries in its conclusion. Since it has no symmetry at all, it also contradicts the weaker statement in which the conclusion asks only for some symmetry other than the identity. 

Since $(\nabla\times B)\times B=(B\cdot\nabla)B-\tfrac12\nabla\|B\|^2$, every solution of \eqref{eq:mhs} is also a steady flow of an ideal incompressible fluid~\cite{Moffatt1985}.

\begin{corollary}\label{cor:euler}
For every $0<\varepsilon\leq1$, the velocity $u=B$ that is tangential to the domain boundary satisfies
\[
 (u\cdot\nabla)u=\nabla\biggl(p+\frac12\|u\|^2\biggr) \, ,\qquad
 \nabla\cdot u=0 \, ,
\]
on \(\Omega_\varepsilon\). The Euler pressure is $-p- (1/2) \|u\|^2$. This flow is real
analytic and nowhere zero, every streamline is closed, and it has no symmetry
other than the identity, in the sense of \eqref{eq:sym} with $u$ in place of $B$. Its Bernoulli function is $-p$, with
critical set $\Gamma$ and nested toroidal level sets.
\end{corollary}

\subsection{Relation to previous work}

Since our equilibria are analytic, their topology is constrained by
Arnold's structure theorem \cite{Arnold65,Arnold1966}; see also
\cite[II, Theorem~1.2]{ArnoldKh}.
In its steady Euler formulation, it states that an analytic flow with
a nonconstant Bernoulli function, defined on a compact analytic
three-dimensional domain and tangent to its boundary, admits,
outside an exceptional analytic set of codimension at least one,
a decomposition into finitely many regions foliated by invariant
tori or invariant cylinders.
The theorem assumes no spatial symmetry. It constrains the topology of
analytic equilibria but does not decide whether nonsymmetric ones
exist, which is the question posed by Conjecture~\ref{conj:grad}. For the flow of Corollary~\ref{cor:euler}, the exceptional set can be
taken to be the axis~$\Gamma$, of codimension two, and
$\Omega_\varepsilon\setminus\Gamma$ is a single region foliated by
invariant tori on which every streamline is closed.



Earlier equilibria without continuous symmetry differ from ours in the solution's symmetry. Lortz~\cite{Lortz1970} proved the existence of toroidal equilibria with nonconstant pressure and no continuous symmetry, in domains symmetric under reflection in a plane and so are consistent with Conjecture~\ref{conj:grad}. Their field lines are closed and the rotational transform is zero. Following Lortz, Drivas, Elgindi and Ginsberg~\cite{DrivasElgindiGinsberg2025} constructed steady Euler flows in a periodic cylinder that retain a plane reflection and a discrete rotation. Salat and Kaiser~\cite{SalatKaiser1995,KaiserSalat1996} found exact equilibria without cylindrical, helical or mirror symmetry in unbounded domains with a straight axis. Bruno and Laurence~\cite{BrunoLaurence1996} and Enciso, Luque and Peralta-Salas~\cite{EncisoLuquePeraltaSalas2025} constructed equilibria without symmetry in toroidal domains with piecewise constant pressure and current sheets. These solutions may be discontinuous across the current sheets.

Recent works give smooth counterexamples to Conjecture~\ref{conj:grad} as stated. Each retains a nontrivial discrete symmetry. Landreman~\cite{Landreman2026} found explicit equilibria in terms of elementary functions with a nowhere-vanishing magnetic field, including a family with rotational transform~$2$ and closed field lines, and a family with sheared rotational transform. His families connect to axisymmetric limits and are stellarator symmetric with two field periods. Issan, Citrin, Landreman and Tracey~\cite{IssanCitrinLandremanTracey2026} removed the stellarator symmetry, but their families retain a half-turn symmetry about the vertical axis. Issan, Jos\'e, Landreman, Citrin and Tracey~\cite{IssanJoseLandremanCitrinTracey2026} found closed-form stellarator-symmetric equilibria with any number of field periods, in which every field line is closed and the rotational transform is an integer. To our knowledge, the equilibria of Theorem~\ref{thm:main} with $\varepsilon>0$ are the first smooth equilibria with nonconstant pressure in a solid torus with nested pressure surfaces whose Euclidean symmetry group is trivial.


Our construction is an expansion about the magnetic axis, in the tradition of Mercier~\cite{Mercier1964}, Solov'ev and Shafranov~\cite{SolovevShafranov1970} and Garren and Boozer~\cite{GarrenBoozer1991a,GarrenBoozer1991b}. Such expansions are formal and are truncated at low order in applications~\cite{LandremanSengupta2018,LandremanSenguptaPlunk2019}. When quasisymmetry is imposed they become overdetermined at third order~\cite{GarrenBoozer1991b}. Weitzner~\cite{Weitzner2014,Weitzner2016} carried formal expansions of nonsymmetric equilibria to all orders in doubly periodic domains, and also expanded about a circular magnetic axis, where his solutions have only finitely many derivatives and neither comes with a proof of convergence. In Theorem~\ref{thm:main} the expansion converges on a neighbourhood of the axis, uniformly in~$\varepsilon$, and its sum is real analytic there, the axis included.


Fixed-boundary codes such as VMEC~\cite{HirshmanWhitson1983} and DESC~\cite{DudtKolemen2020} prescribe the boundary and minimise an energy or a force residual iteratively. In our construction the boundary is instead an output, selected as a pressure level after the series has been summed. Here the data are the circular axis $\Gamma$, the coefficient $X_1$, which fixes the rotational transform $2$, and the harmonic data. Prescribing the boundary in addition to the construction data may overconstrain the problem. More generally, the flux function defining the nested invariant tori and the associated field and pressure profiles must jointly satisfy the scalar normal flux equation \cite{SatoAbe2026}.


\subsection{Ideas of the proof}

Our proof constructs the shape of the pressure surfaces and the pressure profile together. We seek an analytic embedding $X(y,\theta)$, where $y=(y_1,y_2)$ near the origin of $\mathbb R^2$ labels field lines and $\theta\in\mathbb S^1$ runs along them, such that $B\circ X=\partial_\theta X$ and the pressure is independent of $\theta$ in these coordinates. We normalise the pressure to zero on $\Gamma$ and write its pullback as $p(y) = p(X(y,\theta))$ (Section~\ref{sec:coordinates}). The argument indicates which of the two is meant. In these coordinates \eqref{eq:mhs} reduces to a Jacobian condition and two transverse force equations.


The example we pursue is given in Subsection~\ref{subsec:data}. This example is constructed by expanding $X$ about $\Gamma$ and $p$ about zero in homogeneous polynomials of the original transverse labels $y$, and solving one degree at a time. At each degree the transverse frame components of the new coefficient are found by a polynomial Helmholtz--Hodge decomposition, and the tangential component is a primitive in $\theta$. This primitive is periodic only if a mean value vanishes, which is the compatibility condition at the axis. Since every field line is closed, this primitive involves no small denominators. A harmonic vector field may be added freely at each degree. We prescribe a harmonic vector field, of degree two, multiplied by $\varepsilon$. We use an inductive argument to show that this example satisfies the MHD equilibrium equations (Section~\ref{sec:recursion}).


A scalar recursion with a polynomial majorant bounds all coefficients uniformly for $\varepsilon\in[0,1]$ (Sections~\ref{sec:max_frequency}--\ref{sec:majorant}). Consequently, the series converge in the $C^\infty$ topology, that is, uniformly together with all their derivatives on compact subsets of the product of a common label neighbourhood of the axis, the axis included, with a neighbourhood of $[0,1]$ in $\varepsilon$, and their sums are jointly real analytic (Section~\ref{sec:majorant}). Section~\ref{sec:construction} turns the resulting analytic map into an embedding of a neighbourhood of the axis and selects the domain as a pressure superlevel set, using a Morse chart. Subsections~\ref{subsec:pressure-surfaces} and~\ref{subsec:field-lines} then give the geometry of Theorem~\ref{thm:main}. The third derivatives of the pressure on the axis depend on finitely many coefficients and are computed exactly (Appendix~\ref{app:explicit}). For $\varepsilon>0$ they contain the harmonics $\cos3\theta$ and $\sin5\theta$, which exclude every symmetry moving the axis, and a term proportional to $\varepsilon^2$, which excludes the reflection in the horizontal plane (Subsection~\ref{subsec:asymmetry}). The uniform coefficient bounds and the fixed-domain parametrisation give Theorem~\ref{thm:nonisolation} (Subsection~\ref{sec:family}).


The example equilibrium is the sum of a specified convergent series, not an elementary formula or a finite truncation. The figures in Subsection~\ref{sec:figures} show degree-eighteen truncations and are not used in the proof. 


\section{Field-line coordinates and the prescribed data}\label{sec:coordinates}

In this section, we show that an embedding with unit Jacobian
defines an MHD equilibrium if two transverse force balance equations hold,
and we define an example.

\subsection{Reduction of the equilibrium equations}\label{subsec:reduction}

We seek an embedding $X$ of a neighbourhood of $\{0\}\times\mathbb S^1$
in $\mathbb R^2\times\mathbb S^1$ into $\mathbb R^3$.
In transverse labels $y$ and the angle $\theta$, the defining relations are
\begin{equation}\label{eq:BP}
 x = X(y,\theta) \, ,\qquad B( x ) = \partial_\theta X(y,\theta) \, ,
 \qquad p( x ) = p(y) \, .
\end{equation}
Here $y=(y_1,y_2)$ are the transverse labels and $\theta\in\mathbb S^1$.
The notation $p(y)$ denotes the pullback of the physical pressure, which
is independent of $\theta$. Write $\partial_i=\partial/\partial y_i$ for
$i=1,2$ and $\nabla_y=(\partial_1,\partial_2)$.
We use the oriented Jacobian and force functionals
\begin{equation}\label{JFW}
\begin{aligned}
 J&:=\det(\partial_1X,\partial_\theta X,\partial_2X) \, ,\\
 F_i&:=\partial_\theta^2X\cdot\partial_iX
       -\partial_\theta X\cdot\partial_i\partial_\theta X \, ,
       \qquad i=1,2 \, ,\\
 W&:=\partial_1X\cdot\partial_2\partial_\theta X
       -\partial_1\partial_\theta X\cdot\partial_2X \, .
\end{aligned}
\end{equation}

\begin{lemma}\label{lem:reduction}
If $X$ is an embedding with $\partial_\theta J=0$, the equations
$F_i=\partial_i p$, $i=1,2$, imply \eqref{eq:mhs} for the fields in
\eqref{eq:BP}. In particular, it suffices to solve the normalised system
\begin{equation}\label{eq:reduced}
 J=1 \, ,\qquad F_i=\partial_i p \, ,\qquad i=1,2 \, .
\end{equation}
For every $X$, independently of these equations,
\begin{equation}\label{eq:compatibility}
 \partial_2F_1-\partial_1F_2=\partial_\theta W \, .
\end{equation}
\end{lemma}
\begin{proof}
As \(X\) is an embedding, \(J \neq 0\). The field $B$ is the coordinate vector $\partial_\theta$, so its
divergence is $J^{-1}\partial_\theta J$, which vanishes by hypothesis.
The Euclidean identity
$(\nabla\times B)\times B=(B\cdot\nabla)B-\tfrac12\nabla\|B\|^2$
gives
\[
 \bigl(((\nabla\times B)\times B)\circ X\bigr)\cdot\partial_iX
 =\partial_\theta^2X\cdot\partial_iX
       -\tfrac12\partial_i\|\partial_\theta X\|^2=F_i
 =\partial_i p \, .
\]
The last expression equals $(\nabla p\circ X)\cdot\partial_iX$.
Both forces have zero contraction with $\partial_\theta X=B\circ X$:
the magnetic force is perpendicular to $B$, and the pullback of $p$
is independent of $\theta$. Since
$\partial_1X,\partial_\theta X,\partial_2X$ form a basis, the forces agree.
Finally, differentiating \eqref{JFW} and cancelling mixed derivatives gives
\[
 \partial_2F_1-\partial_1F_2
 =\partial_1X\cdot\partial_2\partial_\theta^2X
  -\partial_1\partial_\theta^2X\cdot\partial_2X
 =\partial_\theta W \, .  \tag*{\qedhere}
\]
\end{proof}
Under the normalisation $J=1$, the positive frame and volume order is
$(y_1,\theta,y_2)$, so the argument tuple $(y_1,y_2,\theta)$ is
left-handed; we retain the argument order $X(y,\theta)$ and the displayed
definition of $J$ throughout.

\subsection{Example construction}\label{subsec:data}

Fix an $\varepsilon$. We first construct $X$ in local transverse labels $y$ near the axis. In Subsection~\ref{subsec:unit-disk}, a $\theta$-independent analytic relabelling to the rescaled Morse labels $\xi$ on the unit disk gives the unit-disk embedding $X_{\mathbb D}$ of Theorem~\ref{thm:nonisolation}. All homogeneous coefficients and the normalisation $J=1$ refer to the original labels $y$.

In the labels \((y,\theta)\), we solve \eqref{eq:reduced} first as
formal power series in $y$:
\begin{equation}\label{series}
 X=\sum_{N\geq0}X_N,\qquad p=\sum_{N\geq 0 }p_N.
\end{equation}
Each $X_N=X_N(y,\theta)$ is an $\mathbb R^3$-valued homogeneous
polynomial of degree $N$ in $y$, with $2\pi$-periodic coefficients. Each
$p_N=p_N(y)$ is a real homogeneous polynomial of degree $N$, independent
of $\theta$.
For any scalar polynomial or series $P$, the notation $[P]_N$ extracts
its homogeneous component of degree $N$:
\[
 P=\sum_{j,k\geq0}([y_1^jy_2^k]P)(\theta)y_1^jy_2^k \, ,
 \qquad
 [P]_N=\sum_{j+k=N}([y_1^jy_2^k]P)(\theta)y_1^jy_2^k \, .
\]
Here $[y_1^jy_2^k]P$ denotes a monomial coefficient, and the same
conventions apply componentwise to vectors.
Since $p(y)$ is independent of $\theta$, \eqref{eq:BP} gives $(B\cdot\nabla p)\circ X=\partial_\theta(p\circ X)=0$. This tangency condition also holds for finite truncations of \eqref{series} that are embeddings irrespective of whether they satisfy \eqref{eq:reduced}.

The axis is the unit circle $\Gamma$, parametrised by $X_0=e_1$, where
\[
 e_1=(\cos\theta,\sin\theta,0) \, ,\qquad
 e_\theta=(-\sin\theta,\cos\theta,0) \, ,\qquad e_2=(0,0,1) \, .
\]
This standard cylindrical orthonormal frame has order $(e_1,e_\theta,e_2)$ and satisfies
$e_1\times e_\theta=e_2$,
$\partial_\theta e_1=e_\theta$ and $\partial_\theta e_\theta=-e_1$.
We normalise $p=0$ on $\Gamma$, so that $p_0=0$. Any constant may
subsequently be added to the pressure.
Introduce the linear-form abbreviations
\[
 a_1=y_1\cos2\theta-y_2\sin2\theta,\qquad
 a_2=y_1\sin2\theta+y_2\cos2\theta.
\]
These are not additional coordinates: $y$ remains fixed when
differentiating in $\theta$. They satisfy
$\partial_\theta a_1=-2a_2$, $\partial_\theta a_2=2a_1$ and
$a_1^2+a_2^2=\|y\|^2$. Prescribe
\begin{equation}\label{eq:X1}
 X_1=a_1e_1-a_2e_\theta+a_2e_2 \, .
\end{equation}
At fixed $\theta$, the linear map $y\mapsto(a_1,a_2)$ is a rotation of the label plane through the angle $2\theta$. This double rotation relative to the frame is what produces a rotational transform of $2$. We name the associated differential operators
\[
 \partial_{a_1}=\cos2\theta\,\partial_1-\sin2\theta\,\partial_2,
 \qquad
 \partial_{a_2}=\sin2\theta\,\partial_1+\cos2\theta\,\partial_2,
\]
and write $\nabla_a=(\partial_{a_1},\partial_{a_2})$,
$\nabla_a^\perp=(-\partial_{a_2},\partial_{a_1})$ and, for a planar field
$w=(w_1,w_2)$,
\[
 \nabla_a\cdot w=\partial_{a_1}w_1+\partial_{a_2}w_2,
 \qquad
 \nabla_a\times w=\partial_{a_1}w_2-\partial_{a_2}w_1.
\]
They satisfy $[\partial_\theta,\partial_{a_1}]=-2\partial_{a_2}$ and
$[\partial_\theta,\partial_{a_2}]=2\partial_{a_1}$. At fixed $\theta$ they
have constant coefficients, so they commute with each other and with the
transverse Laplacian
$\Delta_y=\partial_1^2+\partial_2^2=\partial_{a_1}^2+\partial_{a_2}^2$.

For $N\geq2$, we write each polynomial term in coordinates,
\begin{equation}\label{eq:XN}
 X_N=w_1e_1+w_\theta e_\theta+w_2e_2 \, ,
 \qquad w=(w_1,w_2) \, ,
\end{equation}
where $w_1$, $w_\theta$ and $w_2$ are homogeneous polynomials of degree $N$
in $y$ whose coefficients are smooth $2\pi$-periodic functions of $\theta$.
A harmonic datum of degree $N$ is a real periodic homogeneous planar field
$H_N$ of degree $N$ with $\nabla_a\cdot H_N=\nabla_a\times H_N=0$.
Equivalently, $H_N$ is the gradient $\nabla_a$ of the real part of a
complex multiple of $(a_1+\mathrm i a_2)^{N+1}$, with the multiple
depending on $\theta$. With
$g(\theta)=\cos3\theta+\sin5\theta$, we take
\begin{equation}\label{eq:harmonic_datum}
 H_2=\varepsilon\nabla_a\!\left[\frac{g(\theta)}{6}
      \operatorname{Re}\bigl\{(a_1+\mathrm i a_2)^3\bigr\}\right]
     =\frac{\varepsilon g(\theta)}{2}\,(a_1^2-a_2^2,\,-2a_1a_2) \, ,
 \qquad H_N=0\quad(N\geq3) \, .
\end{equation}
The potential is harmonic, so this datum has zero divergence and curl.
The frequencies $3$ and $5$ of $g$ are coprime and occur as a cosine and
a sine. Subsection~\ref{subsec:asymmetry} uses both properties to exclude
every isometry that moves the axis. At leading order, $H_2$ deforms the
circular cross-sections of the pressure surfaces into triangular ones,
with an amplitude that varies along $\Gamma$ as $g$.

For a periodic scalar or polynomial
$P=\sum_{m\in\mathbb Z}P_m(y)e^{\mathrm i m\theta}$, the angular mean and
the zero-mean periodic primitive are
\[
 \langle P\rangle=\frac1{2\pi}\oint_0^{2\pi}P(y,\theta)\,d\theta=P_0 \, ,
 \qquad
 \partial_\theta^{-1}P=\sum_{m\ne0}\frac{P_m(y)e^{\mathrm i m\theta}}{\mathrm i m} \, ,
\]
so that $\partial_\theta\partial_\theta^{-1}P=P-\langle P\rangle$ and
$\langle\partial_\theta^{-1}P\rangle=0$. 

To specify a polynomial right inverse of $\Delta_y$, put
$z=y_1+\mathrm i y_2$ and $\bar z=y_1-\mathrm i y_2$, which are independent
polynomial variables and define $\Delta_y^{-1}$ by a
linear extension of
\begin{equation}\label{eq:laplacian_inverse}
 \Delta_y^{-1}(z^j\bar z^k)
   =\frac{z^{j+1}\bar z^{k+1}}{4(j+1)(k+1)},\qquad j,k\geq0 \, .
\end{equation}
As $\Delta_y=4\partial_z\partial_{\bar z}$, we have
$\Delta_y\Delta_y^{-1}P=P$ for every polynomial $P$. This is a right
inverse only, since $\Delta_y$ annihilates harmonic polynomials. The
operator preserves real polynomials, raises the degree by two and
preserves the angular weight $j-k$. Hence, it commutes with rotations of
the label plane and, at fixed $\theta$, is the same operator when
expressed using the linear forms $a_1,a_2$. It does not commute with
$\partial_1,\partial_2$ on nonzero harmonic polynomials, so the order of
transverse derivatives and $\Delta_y^{-1}$ below matters.

Let $X_{<N}:=X_0+\dots+X_{N-1}$. For $N\geq2$, the degree $N-1$ sources are
\[
 \hat F_i:=[F_i(X_{<N})]_{N-1} \, ,\qquad
 \hat J:=[J(X_{<N})]_{N-1} \, ,\qquad
 \hat W:=[W(X_{<N})]_{N-1} \, ,
\]
the degree-$N$ homogeneous polynomial
\begin{equation}\label{PiNi}
 \mathcal{F}=\frac{y_1\hat F_1+y_2\hat F_2}{N} \, ,
\end{equation}
and
\begin{equation}\label{eq:C0_definition}
 \hat C=\partial_\theta^{-1}
       (\hat W+4\hat J+\partial_{a_1}\mathcal{F}).
\end{equation}
The preliminary normal correction is given by the polynomial
Helmholtz--Hodge decomposition
\begin{equation}\label{eq:normal_preliminary}
 \hat w=-\nabla_a\Delta_y^{-1}\hat J
         +\nabla_a^\perp\Delta_y^{-1}\hat C+H_N
\end{equation}
into a gradient, a co-gradient and the harmonic datum $H_N$, so that
$\nabla_a\cdot\hat w=-\hat J$ and $\nabla_a\times\hat w=\hat C$.

\begin{mydef}[The example]\label{def:example}
Let $X_0=e_1$, let $X_1$ be given by \eqref{eq:X1}, and let $p_0=p_1=0$.
For $N\geq2$, given $X_{<N}$, define successively
\begin{align}
 p_N&=\tfrac23\langle\mathcal{F}-2\hat w_1\rangle \, ,
   \label{eq:pressure_recursion}\\
 w&=\hat w-\tfrac12\nabla_a^\perp\partial_{a_2}\Delta_y^{-1}p_N \, ,
   \label{eq:normal_pressure}\\
 w_\theta&=\partial_\theta^{-1}(\mathcal{F}-p_N-2w_1)  \, ,
   \label{eq:embedding_recursion}
\end{align}
and $X_N$ by \eqref{eq:XN}. The value \eqref{eq:pressure_recursion} is
the one for which the argument of $\partial_\theta^{-1}$ in
\eqref{eq:embedding_recursion} has zero mean, and the last term in
\eqref{eq:normal_pressure} is divergence free. The example is given by the series
\eqref{series} with these coefficients.
\end{mydef}

The coefficients $X_2$, $p_2$ and $p_3$ follow from
Definition~\ref{def:example} by direct computation, and we record them in
Appendix~\ref{app:explicit}. The material from Section~\ref{sec:recursion} to
Subsection~\ref{subsec:embedding} verifies the hypotheses and equations of
Lemma~\ref{lem:reduction} for this example.

\section{Verification of the reduced equations}\label{sec:recursion}

In this section, we prove that the series of Definition~\ref{def:example} solve \eqref{eq:reduced} as formal power series.

\begin{prop}\label{prop:induction}
For every $N$, the coefficient $X_N$ of Definition~\ref{def:example} is a
real periodic homogeneous polynomial of degree $N$, and $p_N$ is a real
$\theta$-independent homogeneous polynomial of degree $N$. Every
homogeneous component of $J-1$ and of $F_i-\partial_ip$, $i=1,2$,
vanishes for the series \eqref{series}.
\end{prop}

\begin{proof}
We expand $F_i$, $J$ and $W$ multilinearly in the homogeneous components
of $X$. In each resulting term, every occurrence of a component of $X$ carries
at most one transverse derivative, and $\partial_i X_0 = 0$, so a nonzero term
containing $X_m$ has degree at least $m-1$. Hence, for $N\geq2$, the components of degree at most $N-2$ depend
only on $X_{<N}$, and those of degree $N-1$ only on $X_{<N}+X_N$.
We say that the induction hypothesis holds at degree $N\geq2$ if $X_{<N}$
is real, \eqref{eq:reduced} holds through degree $N-2$ for $X_{<N}$ and
$p_1+\dots+p_{N-1}$, and the homogeneous components of $W(X_{<N})$ of
degree at most $N-2$ are independent of $\theta$. By the degree count, it
suffices to show that this hypothesis holds at degree $2$ and that, for
each $N\geq2$, it implies the stated properties of $X_N$ and $p_N$ and the
hypothesis at degree $N+1$.

The induction hypothesis holds at degree $2$. Indeed, $X_{<2}=X_0+X_1$ is
real, and the frame components of $\partial_\theta(X_0+X_1)$ and
$\partial_\theta^2(X_0+X_1)$ are $(-a_2,1-a_1,2a_1)$ and
$(-1-a_1,a_2,-4a_2)$. Substituting these and the transverse derivatives
of $X_1$ into \eqref{JFW} gives
\begin{equation}\label{eq:base_functionals}
 F_i(X_0+X_1)=-6(a_1\partial_ia_1+a_2\partial_ia_2)=-6y_i \, ,\quad
 J(X_0+X_1)=1+a_1 \, ,\quad W(X_0+X_1)=-4 \, .
\end{equation}
Hence $[J]_0=1$ and $[F_i]_0=0$. The degree-zero force equations therefore require $\partial_ip_1=0$, so $p_1=0$, and \eqref{eq:reduced} holds in degree zero. Moreover, $[W]_0=-4$ is independent of $\theta$.


Now let $N\geq2$ and suppose that the induction hypothesis holds at degree
$N$. In the components of degree $N-1$ of $F_i(X)$, $J(X)$ and $W(X)$, the
terms without $X_N$ define $\hat F_i$, $\hat J$ and $\hat W$. Since every
occurrence of a component of $X$ carries at most one transverse derivative,
each occurrence of $X_N$ has degree at least $N-1$. A nonzero term
containing $X_N$ at least twice therefore has degree at least $2N-2$, which
exceeds $N-1$ as $N\geq2$. By the same count, in each term of degree $N-1$
containing $X_N$ exactly once, $X_N$ carries a transverse derivative and the
other factors are independent of $y$: they are $\theta$-derivatives of
$X_0$, or transverse derivatives of $X_1$ or of $\partial_\theta X_1$.

In $F_i$ the transverse derivative falls on $X_N$, so the other factor
comes from $X_0$. Using
$\partial_\theta X_N=(\partial_\theta w_1-w_\theta)e_1
+(\partial_\theta w_\theta+w_1)e_\theta+\partial_\theta w_2\,e_2$,
the linear terms are
\[
 \partial_\theta^2X_0\cdot\partial_iX_N
 -\partial_\theta X_0\cdot\partial_i\partial_\theta X_N
 =-e_1\cdot\partial_iX_N-e_\theta\cdot\partial_i\partial_\theta X_N
 =-\partial_i(\partial_\theta w_\theta+2w_1) \, .
\]
Since $\partial_{a_1},\partial_{a_2}$ arise from $\partial_1,\partial_2$ by
a rotation, and $J$ and $W$ are antisymmetric in the two transverse
derivatives,
\[
 J=\det(\partial_{a_1}X,\partial_\theta X,\partial_{a_2}X) \, ,
 \qquad
 W=\partial_{a_1}X\cdot\partial_{a_2}\partial_\theta X
  -\partial_{a_2}X\cdot\partial_{a_1}\partial_\theta X \, .
\]
Moreover, $\partial_{a_1}X_1=e_1$, $\partial_{a_2}X_1=-e_\theta+e_2$,
$\partial_{a_1}\partial_\theta X_1=-e_\theta+2e_2$ and
$\partial_{a_2}\partial_\theta X_1=-e_1$. In $J$ the factor
$\partial_\theta X$ comes from $X_0$, giving $e_\theta$, and the other
transverse derivative falls on $X_1$. As $(e_1,e_\theta,e_2)$ is
positively oriented, the linear terms are
\[
 \det(\partial_{a_1}X_N,e_\theta,-e_\theta+e_2)
 +\det(e_1,e_\theta,\partial_{a_2}X_N)
 =\partial_{a_1}w_1+\partial_{a_2}w_2=\nabla_a\cdot w \, .
\]
In $W$ both factors carry a transverse derivative, so the other factor
comes from $X_1$, and the linear terms are
\[
 e_1\cdot\partial_{a_2}\partial_\theta X_N-e_1\cdot\partial_{a_1}X_N
 +(e_\theta-e_2)\cdot\partial_{a_1}\partial_\theta X_N
 +(e_\theta-2e_2)\cdot\partial_{a_2}X_N \, .
\]
By the frame components of $\partial_\theta X_N$, this equals
$\partial_{a_1}\partial_\theta w_\theta+\partial_{a_2}\partial_\theta w_1
-\partial_{a_1}\partial_\theta w_2-2\partial_{a_2}w_2$, since the terms
$\pm\partial_{a_1}w_1$ and $\pm\partial_{a_2}w_\theta$ cancel. The
commutators give
$-\partial_\theta(\nabla_a\times w)=\partial_{a_2}\partial_\theta w_1
-\partial_{a_1}\partial_\theta w_2+2\nabla_a\cdot w$. Therefore, we have
\begin{equation}\label{XPN}
\begin{aligned}\relax
 [F_i(X)]_{N-1} &= [F_i(X_{<N}+X_N)]_{N-1}
   \hspace{.13em} =\hat F_i-\partial_i(\partial_\theta w_\theta+2w_1),\\
 [J(X)]_{N-1} &= [J(X_{<N}+X_N)]_{N-1}
   \, \, \hspace{.12em} =\hat J+\nabla_a\cdot w,\\
 [W(X)]_{N-1} &= [W(X_{<N}+X_N)]_{N-1}
   =\hat W-\partial_\theta(\nabla_a\times w)-4\nabla_a\cdot w
       +\partial_{a_1}(\partial_\theta w_\theta+2w_1).
\end{aligned}
\end{equation}

The compatibility identity $\partial_2F_1-\partial_1F_2=\partial_\theta W$
holds for every $X$. Its component of degree $N-2$ for $X_{<N}$ and the
induction hypothesis give
$\partial_2\hat F_1-\partial_1\hat F_2=\partial_\theta[W(X_{<N})]_{N-2}=0$.
With Euler's identity $y\cdot\nabla_y\hat F_i=(N-1)\hat F_i$, this gives
\[
 \partial_i(y_1\hat F_1+y_2\hat F_2)
 =\hat F_i+y_1\partial_i\hat F_1+y_2\partial_i\hat F_2
 =\hat F_i+y_1\partial_1\hat F_i+y_2\partial_2\hat F_i=N\hat F_i \, ,
\]
so that $\hat F_i=\partial_i\mathcal{F}$ by \eqref{PiNi}.

We now use \eqref{XPN} to verify that the coefficients of Definition~\ref{def:example} satisfy \eqref{eq:reduced} at degree $N-1$, and that $[W(X)]_{N-1}$ is independent of $\theta$. By
\eqref{eq:normal_pressure},
$w_1=\hat w_1+\tfrac12\partial_{a_2}^2\Delta_y^{-1}p_N$. Since
$\Delta_y^{-1}p_N$ is independent of $\theta$, expanding
$\partial_{a_2}^2$ and using
$\langle\sin^22\theta\rangle=\langle\cos^22\theta\rangle=1/2$ and
$\langle\sin2\theta\cos2\theta\rangle=0$ gives
\[
 \left\langle\partial_{a_2}^2\Delta_y^{-1}p_N\right\rangle
 =\tfrac12\Delta_y\Delta_y^{-1}p_N=\tfrac12p_N \, .
\]
Here the order $\partial_{a_2}\Delta_y^{-1}$ in \eqref{eq:normal_pressure}
is used: the Laplacian acts directly on $\Delta_y^{-1}p_N$. Consequently
$\langle\mathcal{F}-p_N-2w_1\rangle
=\langle\mathcal{F}-2\hat w_1\rangle-\tfrac32p_N=0$ by
\eqref{eq:pressure_recursion}. The argument of $\partial_\theta^{-1}$ in
\eqref{eq:embedding_recursion} therefore has zero mean, so
$\partial_\theta w_\theta+2w_1=\mathcal{F}-p_N$, and the first identity in
\eqref{XPN} becomes
\[
 [F_i(X_{<N}+X_N)]_{N-1}
 =\partial_i\mathcal{F}-\partial_i(\mathcal{F}-p_N)=\partial_i p_N \, .
\]

Since the transverse derivatives commute with each other and $H_N$ has
zero divergence and curl, \eqref{eq:normal_preliminary} gives
$\nabla_a\cdot\hat w=-\Delta_y\Delta_y^{-1}\hat J=-\hat J$ and
$\nabla_a\times\hat w=\Delta_y\Delta_y^{-1}\hat C=\hat C$. The pressure
term in \eqref{eq:normal_pressure} has zero divergence, and its curl is
$-\tfrac12\Delta_y\partial_{a_2}\Delta_y^{-1}p_N=-\tfrac12\partial_{a_2}p_N$
because $\Delta_y$ commutes with $\partial_{a_2}$. Hence
\[
 \nabla_a\cdot w=-\hat J \, ,\qquad
 \nabla_a\times w=\hat C-\tfrac12\partial_{a_2}p_N \, ,
\]
and the second identity in \eqref{XPN} gives
$[J(X_{<N}+X_N)]_{N-1}=0$.

Since $p_N$ is independent of $\theta$, the commutator gives
$\partial_\theta\partial_{a_2}p_N=2\partial_{a_1}p_N$. Together with
\eqref{eq:C0_definition}, this yields
\[
 \partial_\theta(\nabla_a\times w)
 =\hat W+4\hat J+\partial_{a_1}(\mathcal{F}-p_N)
   -\langle\hat W+4\hat J+\partial_{a_1}\mathcal{F}\rangle \, .
\]
Substituting this, $\nabla_a\cdot w=-\hat J$ and
$\partial_\theta w_\theta+2w_1=\mathcal{F}-p_N$ into the third identity in
\eqref{XPN} gives
$[W(X_{<N}+X_N)]_{N-1}=\langle\hat W+4\hat J+\partial_{a_1}\mathcal{F}\rangle$,
which is independent of $\theta$.

Finally, $\hat J$, $\hat W$ and $\hat C$ have degree $N-1$ and
$\mathcal{F}$ has degree $N$. Since $\Delta_y^{-1}$ raises the degree by
two and each transverse derivative lowers it by one, $\hat w$, $w$,
$w_\theta$ and $p_N$ are homogeneous of degree $N$. Since $X_{<N}$ and
$H_N$ are real and every operation used preserves reality and
periodicity, $X_N$ and $p_N$ are real and periodic, and averaging makes
$p_N$ independent of $\theta$. As $\partial_ip_N$ has degree $N-1$ and
adding $X_N$ changes no component of degree at most $N-2$, the induction
hypothesis holds at degree $N+1$.
\end{proof}

\begin{remark}[Non-uniqueness]
Besides the choice of $H_N$, one may add to $w$ any field
$\nabla_a^\perp\Delta_y^{-1}P$, with $P$ a real $\theta$-independent
homogeneous polynomial of degree $N-1$, and recompute $w_\theta$ from
\eqref{eq:embedding_recursion}. This leaves $p_N$ and the conclusions of
Proposition~\ref{prop:induction} unchanged. At degree $N$ this addition is
the effect of a $\theta$-independent area-preserving relabelling of the
transverse labels, and our construction omits it.
\end{remark}

\section{Linear bound on the maximum \texorpdfstring{$\theta$}{theta}-frequency}\label{sec:max_frequency}

In this section, we show that the largest angular Fourier frequency in
$X_N$ grows at most linearly with $N$.
Section~\ref{sec:bound_degree} uses this bound to control second angular
derivatives, and Section~\ref{sec:majorant} uses it to show that $X$ is
real analytic in $\theta$, and not merely smooth, as required for the
joint analyticity asserted in Theorem~\ref{thm:nonisolation}.

\begin{prop}[Maximum frequency]\label{prop:frequency}
Let $f_N$ be the largest absolute $\theta$-frequency in the Cartesian
components of $X_N$ given via Definition \ref{def:example}. Then
\[
 f_0=1,\qquad f_1 \leq 3,\qquad f_2\leq14,\qquad
 f_N\leq32N-48\quad(N\geq2),\qquad f_N\leq32N\quad(N\geq1).
\]
\end{prop}
\begin{proof}
The proof has three steps. We first record how the operations in
Definition~\ref{def:example} change frequencies, and deduce that the solve
from the data to $X_N$ raises the maximum frequency by at most
five. We then compute $f_0$, $f_1$ and $f_2$ directly. Finally, we prove $f_N\leq32N-48$ by induction for $N\geq2$,
with $f_2$ as the base case.

The maximum frequency of a product is at most the sum of the maxima
of its factors, and that of a sum is at most their maximum.
Differentiation in $\theta,y_1,y_2$ does not increase frequency.
Neither do $\Delta_y^{-1}$, averaging, or $\partial_\theta^{-1}$:
they act separately on Fourier modes, and averaging keeps only the
zero mode. Each $\partial_{a_i}$ can increase frequency by at most two.
Since $X_0$, $X_1$ and the data \eqref{eq:harmonic_datum} are Fourier
polynomials, so is every $X_N$.

The solve from $\mathcal{F},\hat J,\hat W,H_N$ to $X_N$ increases their
maximum frequency by at most five. To see this, let $\hat f_N$ denote the
largest absolute $\theta$-frequency in these data.
Equation~\eqref{eq:C0_definition} gives at most $\hat f_N+2$ for $\hat C$,
and \eqref{eq:normal_preliminary} gives at most $\hat f_N+4$ for $\hat w$.
The pressure $p_N$ and its potential $\Delta_y^{-1}p_N$ are independent
of $\theta$, so the two $a$-derivatives in the pressure correction
give frequency at most four. Hence $w$ has frequency at most $\hat f_N+4$.
The primitive in \eqref{eq:embedding_recursion} gives the same bound
for $w_\theta$. Finally, multiplication by $e_1$ or $e_\theta$ adds
at most one and $e_2$ is constant. This proves $f_N\leq\hat f_N+5$.

We have $f_0=1$, and $f_1\leq3$ since the frame components of $X_1$ in
\eqref{eq:X1} have frequency two. By \eqref{eq:base_functionals} and
\eqref{PiNi}, at degree two the sources $\mathcal{F}=-3\|y\|^2$,
$\hat J=a_1$ and $\hat W=0$ have frequencies at most two. The components
of $H_2$ in \eqref{eq:harmonic_datum} have frequency at most nine: $g$ has
maximum frequency five, and $a_1^2-a_2^2$ and $a_1a_2$ have frequency
four. Hence $\hat f_2\leq9$, and the preceding solve gives
$f_2\leq\hat f_2+5\leq14$.

We prove $f_N\leq32N-48$ by induction for $N\geq2$. The offset $-48$
strengthens the induction hypothesis. A source term may have total index
$N+1$, which exceeds $N$, so the weaker bound $f_j\leq32j$ alone would
not close the induction. The saving of $48$ in each factor of index at
least two compensates for this excess and for the increase of five in
the solve. The base case holds because $14\leq32\cdot2-48$. Suppose that
$f_j\leq32j-48$ for $2\leq j\leq N-1$, where $N\geq3$. Since $H_N=0$ for
$N\geq3$ by \eqref{eq:harmonic_datum}, only the known sources need to be
counted.
In $\mathcal{F}$, every term has two factors with degree indices $j+k=N$:
the transverse derivative in $\hat F_i$ is offset by multiplication
by $y_i$. In $\hat W$, two transverse derivatives give $j+k=N+1$.
In $\hat J$, terms have the form
\[
 \det(\partial_1X_j,\partial_\theta X_k,\partial_2X_l),
 \qquad j+k+l=N+1.
\]
All indices are less than $N$, since the sources use $X_{<N}$.
Transversely differentiated factors have positive indices because
$X_0$ is independent of $y$.

Each nonzero source term therefore has at most three factors, total
index at most $N+1$, and at most one zero index. At least one index
is at least two: otherwise the index sums would be at most two for
$\mathcal{F},\hat W$, and at most three for $\hat J$, contradicting their
respective constraints when $N\geq3$.
If exactly one factor has index $j\geq2$, the other factors have
frequencies at most $\max\{f_0,f_1\}\leq3$, and $j\leq N-1$.
Such a term therefore has frequency at most
\[
 f_j+3+3\leq32(N-1)-48+6\leq32N-53.
\]
If at least two indices, say $j,k$, are at least two, their factors
save $48$ each relative to the weaker bound $f_j\leq32j$:
$f_j+f_k\leq32(j+k)-96$.
Every other positive index satisfies the weaker bound as well, and a
possible zero index contributes at most one. As the total index is
at most $N+1$, such a term has frequency at most
\[
 32(N+1)-2\cdot48+1\leq32N-53.
\]
Taking the maximum over all source terms gives $\hat f_N\leq32N-53$,
and hence $f_N\leq\hat f_N+5\leq32N-48$. This completes the induction
and proves $f_N\leq32N$ for $N\geq1$.
\end{proof}

\section{Bound on the degree \texorpdfstring{$N$}{N} correction}\label{sec:bound_degree}

In this section, we derive a scalar inequality that bounds each
correction in terms of the prescribed harmonic datum and the
corrections of lower degree.
In Section~\ref{sec:majorant}, this inequality gives a radius of
convergence independent of $\varepsilon$, as needed for the proof
of Theorem~\ref{thm:nonisolation}. 

\subsection{Polynomial norm}

In this subsection, we establish estimates for products, derivatives
and the right inverse of $\Delta_y$ in the norms used to measure the
homogeneous coefficients.

For a homogeneous polynomial $P_N$ of degree $N$ in the independent
variables $z,\bar z$, define
\begin{equation}\label{eq:polynomial-norm}
 \|P_N\|_*=
 \max_{0\leq n\leq N}
 \frac{ \big|[z^n\bar z^{N-n}]P_N \big|}{\binom Nn} \, .
\end{equation}
Here $[z^n\bar z^{N-n}]$ extracts the indicated monomial coefficient.

\begin{prop}[Product norm]
For homogeneous polynomials $P_j,P_k$ of degrees $j,k$,
\[
 \|P_jP_k\|_*\leq\|P_j\|_*\|P_k\|_*.
\]
\end{prop}
\begin{proof}
The definition gives
$|[z^l\bar z^{j-l}]P_j|\leq\binom jl\|P_j\|_*$ and
$|[z^{n-l}\bar z^{k-n+l}]P_k|
 \leq\binom k{n-l}\|P_k\|_*$.
Consequently the coefficient of $z^n\bar z^{j+k-n}$ in their product satisfies
\[
\begin{aligned}
 \big| [z^n\bar z^{j+k-n}](P_jP_k) \big|
 &\leq\sum_l
 \big|[z^l\bar z^{j-l}]P_j \big| \,
 \big| [z^{n-l}\bar z^{k-n+l}]P_k \big|\\
 &\leq\|P_j\|_*\|P_k\|_*
 \sum_l\binom jl\binom k{n-l}
 =\binom{j+k}{n}\|P_j\|_*\|P_k\|_*.
\end{aligned}
\]
The sum ranges over $0\leq l\leq j$ and $0\leq n-l\leq k$. The last equality is Vandermonde's identity. Dividing by
$\binom{j+k}{n}$ and taking the maximum over $n$ proves the claim.
\end{proof}

\begin{prop}[Transverse derivatives and the right inverse]\label{prop:transversal}
For a homogeneous polynomial $P_N$,
\[
 \|\partial_iP_N\|_*\leq2N\|P_N\|_*\quad(i=1,2),
 \qquad
 \|\Delta_y^{-1}P_N\|_*
 \leq\frac{\|P_N\|_*}{4(N+1)(N+2)} \, .
\]
\end{prop}
\begin{proof}
For $N\geq1$, coefficient extraction gives
\[
 \|\partial_zP_N\|_*
 =\max_{1\leq n\leq N}
 \frac{n \big|[z^n\bar z^{N-n}]P_N \big| }{\binom{N-1}{n-1}}
 =N\max_{1\leq n\leq N}
 \frac{ \big| [z^n\bar z^{N-n}]P_N \big| }{\binom Nn}
 \leq N\|P_N\|_* \, .
\]
The same argument gives $\|\partial_{\bar z}P_N\|_*\leq N\|P_N\|_*$.
Since $\partial_1=\partial_z+\partial_{\bar z}$ and
$\partial_2=\mathrm i(\partial_z-\partial_{\bar z})$, the triangle
inequality gives both transverse derivative bounds. For $N=0$ the
derivatives vanish.

The definition of the right inverse gives
\[
 \Delta_y^{-1}P_N
 =\sum_{n=0}^N
 \frac{ \big( [z^n\bar z^{N-n}]P_N \big) z^{n+1}\bar z^{N-n+1}}
 {4(n+1)(N-n+1)} \, .
\]
This has degree $N+2$, so the binomial weight of its displayed
coefficient is $\binom{N+2}{n+1}$. The identity
\[
 \binom{N+2}{n+1}
 =\frac{(N+2)!}{(n+1)!(N-n+1)!}
 =\frac{(N+1)(N+2)}{(n+1)(N-n+1)}\binom Nn
\]
therefore implies
\[
 \|\Delta_y^{-1}P_N\|_*
 =\max_{0\leq n\leq N}
 \frac{ \big| [z^n\bar z^{N-n}]P_N \big| }
 {4(n+1)(N-n+1)\binom{N+2}{n+1}}
 =\frac{\|P_N\|_*}{4(N+1)(N+2)} \, . \tag*{\qedhere}
\]
\end{proof}

For a Fourier polynomial $P(y,\theta)=\sum_mP_m(y)e^{\mathrm im\theta}$
whose coefficients have one common homogeneous degree, define
\begin{equation}\label{eq:fourier-norm}
 \|P\|_F=\sum_m\|P_m\|_* \, .
\end{equation}
The combination that measures one angular derivative is simply
\begin{equation}\label{eq:one-derivative-norm}
 \|P\|_F+\|\partial_\theta P\|_F
 =\sum_m \big(1+|m| \big) \|P_m\|_* \, .
\end{equation}
For vector polynomials these coefficient norms are summed over
components. For $X_N$ the components are in a fixed Cartesian basis.

\begin{prop}[Fourier norm estimates]\label{prop:Fourier}
Let $P$ be a Fourier polynomial with frequencies $|m|\leq f$.
Then
\[
\begin{aligned}
 \|\partial_\theta P\|_F
 &\leq\|P\|_F+\|\partial_\theta P\|_F,\\
 \|\partial_\theta^2P\|_F
 &\leq f(\|P\|_F+\|\partial_\theta P\|_F),\\
 \|\partial_\theta P\|_F+\|\partial_\theta^2P\|_F
 &\leq f(\|P\|_F+\|\partial_\theta P\|_F),\\
 \|\partial_\theta^{-1}P\|_F+
 \|\partial_\theta\partial_\theta^{-1}P\|_F
 &\leq2\|P\|_F,\qquad
 \|\langle P\rangle\|_F\leq\|P\|_F \, .
\end{aligned}
\]
When applied to $X_N$ with $N\geq1$, Proposition~\ref{prop:frequency} supplies $f=f_N\leq32N$.
\end{prop}
\begin{proof}
The Fourier coefficients of $\partial_\theta P$ and
$\partial_\theta^2P$ are $\mathrm imP_m$ and $-m^2P_m$.
The first three estimates follow by summing
\[
 |m|\leq1+|m|,\qquad m^2\leq f(1+|m|),\qquad
 |m|(1+|m|)\leq f(1+|m|).
\]
Moreover,
$\partial_\theta^{-1}P=\sum_{m\ne0}P_m e^{\mathrm im\theta}/(\mathrm im)$
and $\langle P\rangle=P_0$. Hence
\begin{align*}
 \|\partial_\theta^{-1}P\|_F+
 \|\partial_\theta\partial_\theta^{-1}P\|_F
 &=\sum_{m\ne0}\frac{1+|m|}{|m|}\|P_m\|_*
 \leq2\sum_{m\ne0}\|P_m\|_*\leq2\|P\|_F,\\
 \|\langle P\rangle\|_F&=\|P_0\|_*\leq\sum_m\|P_m\|_*=\|P\|_F. \tag*{\qedhere}
\end{align*}
\end{proof}

\subsection{Bound on the degree \texorpdfstring{$N$}{N} correction}


In this subsection, we measure the degree-$N$ correction $(X_N,p_N)$ in
the norms of the previous subsection and bound it in terms of the harmonic
datum and the corrections of lower degree.
In Section~\ref{sec:majorant} we utilise this for our convergence proof.



Define the sequence
\begin{equation}\label{eq:size-sequence}
 \beta_N=N\bigl(\|X_N\|_F+\|\partial_\theta X_N\|_F+\|p_N\|_F\bigr)
 \qquad(N\geq1).
\end{equation}

\begin{prop}[Bound on the degree $N$ correction]\label{prop:bound}
For $N\geq2$, the data $\mathcal{F},\hat J,\hat W,H_N$
determine a correction satisfying
\begin{equation}\label{eq:bound}
 \frac{\beta_N}{N}\leq C\left[
 \|\mathcal{F}\|_F+
 \frac{\|\hat J\|_F+\|\partial_\theta\hat J\|_F+\|\hat W\|_F}{N}
 +\|H_N\|_F+\|\partial_\theta H_N\|_F\right],
\end{equation}
where $C\geq1$ is independent of $N$ and the data.
\end{prop}
\begin{proof}
In this proof $C$ may be enlarged between estimates. 
Proposition~\ref{prop:transversal} applies mode by mode. For a
degree-$k$ polynomial $P_k$ and $i_1,i_2\in\{1,2\}$ it gives
\[
\begin{aligned}
 \|\partial_{i_1}\Delta_y^{-1}P_k\|_F
 &\leq2(k+2)\|\Delta_y^{-1}P_k\|_F
 \leq\frac{\|P_k\|_F}{2(k+1)},\\
 \|\partial_{i_1}\partial_{i_2}\Delta_y^{-1}P_k\|_F
 &\leq4(k+1)(k+2)\|\Delta_y^{-1}P_k\|_F
 \leq\|P_k\|_F.
\end{aligned}
\]
The same bounds hold after applying $\partial_\theta$, since the
displayed operators commute with it. For $\partial_{a_i}$,
multiplication by $\cos2\theta$ or $\sin2\theta$ does not increase
$\|\cdot\|_F$, and the commutators
$[\partial_\theta,\partial_{a_1}]=-2\partial_{a_2}$ and
$[\partial_\theta,\partial_{a_2}]=2\partial_{a_1}$ add only terms of
the same form, so these estimates change by an absolute constant,
no commutation of
$\Delta_y^{-1}$ with transverse differentiation is used.

The definition
$\hat C=\partial_\theta^{-1}
(\hat W+4\hat J+\partial_{a_1}\mathcal{F})$ and the primitive estimate yield
\[
 \|\hat C\|_F+\|\partial_\theta\hat C\|_F
 \leq C\bigl(\|\hat W\|_F+\|\hat J\|_F+N\|\mathcal{F}\|_F\bigr).
\]
Both $\hat J$ and $\hat C$ have degree $N-1$. Applying the first
transverse estimate to the Hodge formula~\eqref{eq:normal_preliminary},
\[
\begin{aligned}
 \|\hat w\|_F+\|\partial_\theta\hat w\|_F
 &\leq\frac{C}{N}
 \bigl(\|\hat J\|_F+\|\partial_\theta\hat J\|_F+
       \|\hat C\|_F+\|\partial_\theta\hat C\|_F\bigr)
 +C(\|H_N\|_F+\|\partial_\theta H_N\|_F)\\
 &\leq C\left[\|\mathcal{F}\|_F+
 \frac{\|\hat J\|_F+\|\partial_\theta\hat J\|_F+\|\hat W\|_F}{N}
 +\|H_N\|_F+\|\partial_\theta H_N\|_F\right].
\end{aligned}
\]
The averaging formula~\eqref{eq:pressure_recursion} gives
\[
 \|p_N\|_F\leq C(\|\mathcal{F}\|_F+\|\hat w\|_F).
\]
The pressure correction in~\eqref{eq:normal_pressure} is
$-\tfrac12\nabla_a^\perp\partial_{a_2}\Delta_y^{-1}p_N$.
Its components involve two transverse derivatives of the right inverse
and fixed trigonometric factors. Since $p_N$ is independent of $\theta$,
the second transverse estimate gives
\[
 \|w-\hat w\|_F+\|\partial_\theta(w-\hat w)\|_F
 \leq C\|p_N\|_F,
 \qquad w=(w_1,w_2).
\]
Finally,~\eqref{eq:embedding_recursion} and the primitive estimate imply
\[
\begin{aligned}
 \|w_\theta\|_F+\|\partial_\theta w_\theta\|_F
 &\leq C(\|\mathcal{F}\|_F+\|p_N\|_F+\|w\|_F),\\
 \|X_N\|_F+\|\partial_\theta X_N\|_F
 &\leq C\bigl(\|w_\theta\|_F+\|\partial_\theta w_\theta\|_F
              +\|w\|_F+\|\partial_\theta w\|_F\bigr).
\end{aligned}
\]
The last line uses the Cartesian component convention and the bounded
derivatives of the frame. Combining these estimates proves the claim.
\end{proof}

To cover a complex neighbourhood of $[0,1]$, let $C_H$ be four times the value of
$C(\|H_2\|_F+\|\partial_\theta H_2\|_F)$ at $\varepsilon=1$,
where $C$ is the constant in~\eqref{eq:bound}, which does not depend on
$C_H$. By~\eqref{eq:harmonic_datum}, $H_2$ is linear in $\varepsilon$ and
$H_N=0$ for $N\geq3$, so after multiplying~\eqref{eq:bound} by $N$ the
harmonic datum contributes at most $C_H\delta_{N,2}$ for
$|\varepsilon|<2$, where $\delta_{N,2}$ is the Kronecker delta.

\newcommand{\lnormlow}[2]{\mathopen{\raisebox{#2}{$\nulldelimiterspace=0pt
  \left\lVert\vcenter{\hbox{\rule{0pt}{#1}}}\right.$}}}
\newcommand{\rnormlow}[2]{\mathclose{\raisebox{#2}{$\nulldelimiterspace=0pt
  \left.\vcenter{\hbox{\rule{0pt}{#1}}}\right\rVert$}}}

\begin{prop}[Bound via lower-degree corrections]
For every $N\geq2$ and $|\varepsilon|<2$,
\begin{equation}\label{eq:beta}
 \beta_N\leq C_H\delta_{N,2}
 +C\left[
 \sum_{j+k=N}\beta_j\beta_k
 +\sum_{\substack{j+k=N+1\\j,k\geq2}}\beta_j\beta_k
 +\sum_{j+k+l=N+1}\beta_j\beta_k\beta_l\right].
\end{equation}
Unless stated otherwise, the summation indices are positive integers.
\end{prop}
\begin{proof}
We estimate each known source in~\eqref{eq:bound}. Homogeneity gives
\[
\begin{aligned}
 \|\mathcal{F}\|_F
 &=\frac1N \, \|y_1\hat F_1+y_2\hat F_2\|_F\\
 &=
 \frac1N \, \lnormlow{4em}{-1.25ex}
 \sum_{\substack{j+k=N\\j,k\geq1}}
 \sum_{i=1}^2 y_i
 \bigl(\partial_\theta^2X_j\cdot\partial_iX_k
       -\partial_\theta X_j\cdot\partial_i\partial_\theta X_k\bigr)
 \rnormlow{4em}{-1.25ex}_F \\ 
 &=\frac1N  \left\|
 \sum_{j+k=N}k\bigl(\partial_\theta^2X_j\cdot X_k
                 -\partial_\theta X_j\cdot\partial_\theta X_k\bigr)
 \right\|_F\\
 &\leq\frac C N\sum_{j+k=N}k(1+32j)
                      \frac{\beta_j}{j}\frac{\beta_k}{k}.
\end{aligned}
\]
All factors come from $X_{<N}$ and a zero index cannot contribute since
the other index would be $N$. The third line uses
$y\cdot\nabla_yX_k=kX_k$, and the last uses
Proposition~\ref{prop:Fourier} with $f_j\leq32j$ from
Proposition~\ref{prop:frequency}.

For the Jacobian source, the orientation convention gives exactly
\[
 \hat J=
 \sum_{j,k,l=0}^{N-1}
 [\det(\partial_1X_j,\partial_\theta X_k,\partial_2X_l)]_{N-1} \, .
\]
The terms with $j=0$ or $l=0$ vanish. Separating $k=0$ from $k\geq1$
and applying the product bounds, including one angular derivative, gives
\[
\begin{aligned}
 \|\hat J\|_F+\|\partial_\theta\hat J\|_F
 &\leq C\sum_{\substack{j+l=N+1\\j,l\geq2}}
 jl\frac{\beta_j}{j}\frac{\beta_l}{l}
 +32C\sum_{j+k+l=N+1}
 jkl\frac{\beta_j}{j}\frac{\beta_k}{k}\frac{\beta_l}{l}\\
 &=C\sum_{\substack{j+l=N+1\\j,l\geq2}}\beta_j\beta_l
 +32C\sum_{j+k+l=N+1}\beta_j\beta_k\beta_l.
\end{aligned}
\]
In the first sum $j,l\leq N-1$ and $j+l=N+1$ force $j,l\geq2$.
The fixed norms of $\partial_\theta X_0$ and
$\partial_\theta^2X_0$ are absorbed into $C$. In the second sum,
Proposition~\ref{prop:Fourier} bounds
$\|\partial_\theta X_k\|_F+\|\partial_\theta^2X_k\|_F$
by $32\beta_k$, and each transverse derivative contributes its
degree factor. Any extra terms included by dropping the upper index
bounds are nonnegative.

Similarly, relabelling the two factors in the second term of $W$ gives
\[
\begin{aligned}
 \|\hat W\|_F
 &=\left\|
 \sum_{j,k=0}^{N-1}
 [\partial_1X_j\cdot\partial_2\partial_\theta X_k
 -\partial_2X_j\cdot\partial_1\partial_\theta X_k]_{N-1}
 \right\|_F\\
 &=\lnormlow{3.6em}{-1.5ex}
 \sum_{\substack{j+k=N+1\\j,k\geq2}}
 (\partial_1X_j\cdot\partial_2\partial_\theta X_k
 -\partial_2X_j\cdot\partial_1\partial_\theta X_k)
 \rnormlow{3.6em}{-1.5ex}_F \\
 &\leq C\sum_{\substack{j+k=N+1\\j,k\geq2}}
 jk\frac{\beta_j}{j}\frac{\beta_k}{k}
 =C\sum_{\substack{j+k=N+1\\j,k\geq2}}\beta_j\beta_k.
\end{aligned}
\]
Both factors are transversely differentiated, so positive indices are
required, and the degree constraint and truncation strengthen them to
$j,k\geq2$. Only the already controlled first angular derivative is
used here.

Substitute these bounds into~\eqref{eq:bound} and multiply by $N$.
After interchanging $j,k$ in the force-source sum, its index factor is
removed by
\[
 j(1+32k)\frac{\beta_j}{j}\frac{\beta_k}{k}
 =\frac{\beta_j\beta_k}{k}+32\beta_j\beta_k
 \leq33\beta_j\beta_k.
\]
The other two sums already have no index weights. Enlarging $C$ to
absorb the constants $33$ and $32$ proves~\eqref{eq:beta}.
\end{proof}

\section{A convergent scalar series dominating every coefficient}\label{sec:majorant}

In this section, we use a scalar majorant to show that the series for
$X$ and $p$ converge to jointly analytic functions on a neighbourhood
of $\{0\}\times\mathbb S^1\times[0,1]$ in $(y,\theta,\varepsilon)$.
This common neighbourhood allows us to construct the family of analytic
embeddings in Section~\ref{sec:construction} and prove
Theorem~\ref{thm:nonisolation} in Subsection~\ref{sec:family}.

Every coefficient is polynomial in $\varepsilon$: the datum
\eqref{eq:harmonic_datum} is linear in it, products preserve polynomial dependence, and all linear
operators are independent of it. The
norm estimates hold also for complex coefficients. We work on
$|\varepsilon|<2$, a complex neighbourhood of $[0,1]$, and use the
uniform forcing bound already incorporated in $C_H$.

In every sum in~\eqref{eq:beta} all indices lie in
$\{1,\ldots,N-1\}$. This follows from $j+k=N$ with $j,k\geq1$,
from $j+k=N+1$ with $j,k\geq2$, and from $j+k+l=N+1$ with
$j,k,l\geq1$, respectively. Since $X_1$ does not depend on
$\varepsilon$, neither does $\beta_1$. Define $b_1=\beta_1$ and
\begin{equation}\label{eq:w}
 b_N=C_H\delta_{N,2}
 +C\left[
 \sum_{j+k=N}b_jb_k
 +\sum_{\substack{j+k=N+1\\j,k\geq2}}b_jb_k
 +\sum_{j+k+l=N+1}b_jb_kb_l\right]
 \qquad(N\geq2).
\end{equation}
Every $b_N$ is nonnegative. The right-hand side is a polynomial with
nonnegative coefficients in the preceding terms, and is therefore
nondecreasing in each of them on $[0,\infty)$.

Induction now gives
\begin{equation}\label{eq:a-w}
 \beta_N\leq b_N\qquad(N\geq1).
\end{equation}
For $N=1$ this is equality. For $N\geq2$, every index in
\eqref{eq:beta} is less than $N$, so the induction hypothesis,
positivity and \eqref{eq:w} give the bound directly.
From this point until the geometric bound is established, only the
nonnegative scalar recursion~\eqref{eq:w} is used.

Let $t$ be a bookkeeping variable. For a
trial polynomial $q(t)=\sum_{n=1}^dq_nt^n$ of finite degree with $q_n\geq0$ and $q_1=b_1$, define
\[
 \mathcal T[q](t)=C_Ht^2+Cq(t)^2
 +Ct^{-1}(q(t)-b_1t)^2+Ct^{-1}q(t)^3.
\]
Because $t^2$ divides $q(t)-b_1t$ and $t^3$ divides $q(t)^3$,
this is a polynomial with nonnegative coefficients and no terms of
degree zero or one. Set $q_n=0$ beyond the degree of $q$, and let
$[t^n]q$ denote coefficient extraction. The identities
\[
 [t^n]q^2=\sum_{j+k=n}q_jq_k,\qquad \quad
 [t^n]t^{-1}(q-b_1t)^2
 =[t^{n+1}]\left(\sum_{j\geq2}q_jt^j\right)^2,
 \qquad \quad [t^n]t^{-1}q^3=[t^{n+1}]q^3
\]
give, for $n\geq2$,
\begin{equation}\label{eq:G-coeff}
 [t^n]\mathcal T[q]=C_H\delta_{n,2}
 +C\left[
 \sum_{j+k=n}q_jq_k
 +\sum_{\substack{j+k=n+1\\j,k\geq2}}q_jq_k
 +\sum_{j+k+l=n+1}q_jq_kq_l\right].
\end{equation}
The same index count shows that this depends only on
$q_1,\ldots,q_{n-1}$. At a positive evaluation radius $R_*$,
the value of $\mathcal T[q]$ depends on $q$ only through $q(R_*)$:
\begin{equation}\label{eq:G-value}
 \mathcal T[q](R_*)=C_HR_*^2+Cq(R_*)^2
 +CR_*^{-1}(q(R_*)-b_1R_*)^2+CR_*^{-1}q(R_*)^3.
\end{equation}
Define the partial generating polynomials
$\mathcal S_N(t)=\sum_{n=1}^Nb_nt^n$. For $2\leq n\leq N$, the
first $n-1$ coefficients of $\mathcal S_{N-1}$ are
$b_1,\ldots,b_{n-1}$. Comparing~\eqref{eq:G-coeff} with~\eqref{eq:w}
therefore gives
\begin{equation}\label{eq:G-reproduces}
 [t^n]\mathcal T[\mathcal S_{N-1}]=b_n
 \qquad(2\leq n\leq N).
\end{equation}
We apply $\mathcal T$ only to the polynomials $\mathcal S_{N-1}$, so
trial polynomials of finite degree suffice. Working with polynomials
rather than the full generating series avoids assuming in advance that
$\sum_{n\geq1}b_nR_*^n$ converges.

\begin{lemma}\label{lem:majorant}
If $\mu>0$ and $R_*>0$ satisfy
\begin{equation}\label{eq:eps-R}
 C_HR_*+CR_*(b_1+\mu)^2+C\mu^2
 +CR_*(b_1+\mu)^3\leq\mu,
\end{equation}
then $\sum_{n\geq1}b_nR_*^n\leq(b_1+\mu)R_*$.
\end{lemma}
\begin{proof}
We prove $\mathcal S_N(R_*)\leq(b_1+\mu)R_*$ by induction.
For $N=1$, this is $b_1R_*\leq(b_1+\mu)R_*$. For $N\geq2$,
\eqref{eq:G-reproduces} and nonnegativity of all coefficients give
\[
 \mathcal S_N(R_*)
 =b_1R_*+\sum_{n=2}^N
 ([t^n]\mathcal T[\mathcal S_{N-1}])R_*^n
 \leq b_1R_*+\mathcal T[\mathcal S_{N-1}](R_*).
\]
Write $\mathcal S_{N-1}(R_*)=(b_1+\kappa)R_*$. Nonnegativity gives
$\kappa\geq0$, and the induction hypothesis gives $\kappa\leq\mu$.
Equation~\eqref{eq:G-value} yields
\[
 \mathcal T[\mathcal S_{N-1}](R_*)
 =R_*\left[C_HR_*+CR_*(b_1+\kappa)^2+C\kappa^2
              +CR_*(b_1+\kappa)^3\right]\leq\mu R_*.
\]
The bracket is nondecreasing for $\kappa\geq0$, and the last
inequality follows from~\eqref{eq:eps-R}. Thus
$\mathcal S_N(R_*)\leq(b_1+\mu)R_*$. The partial sums of the
nonnegative series are bounded by this same quantity, and hence so
is its sum.
\end{proof}
Choose
\[
 \mu=\frac1{2C},\qquad
 R_*=
 \frac{\mu}{2 \big( C_H+C(b_1+\mu)^2+C(b_1+\mu)^3 \big)}.
\]
Then $C\mu^2=\mu/2$, while the other three terms on the left of
\eqref{eq:eps-R} sum to
$R_*[C_H+C(b_1+\mu)^2+C(b_1+\mu)^3]=\mu/2$.
Thus~\eqref{eq:eps-R} holds with equality. Each term of a convergent
nonnegative series is at most its sum. Lemma~\ref{lem:majorant}
therefore gives $b_NR_*^N\leq(b_1+\mu)R_*$. Together with
\eqref{eq:a-w}, this proves
\begin{equation}\label{eq:uniform-majorant}
 \beta_N\leq b_N\leq MR_*^{-N},
 \qquad M=(b_1+\mu)R_*,\qquad N\geq1.
\end{equation}
Both $M$ and $R_*$ are independent of $\varepsilon$. The latter is
a comparison-series radius, distinct from the physical cylindrical
radius $R$. Its formula is explicit in $C,b_1,C_H$, but does not
give a numerical radius because the constant $C$ has not been computed.

We next prove joint analyticity. The construction and all estimates
so far use $\theta\in\mathbb S^1$. Here we extend to a complex $\theta$ to establish analyticity in $\theta$. Each $X_N$ is a finite
Fourier sum, so it extends to $\theta\in\mathbb C/2\pi\mathbb Z$ by the
same formula, with $\mathbb S^1$ the real slice
$\operatorname{Im}\theta=0$. We utilise: a function
holomorphic on an open set containing the real points restricts to a
real-analytic function of the real variables. Complexify $y$ so that
$z,\bar z$ are independent complex variables, and extend $\theta$ in
this way.

Each Fourier factor satisfies
$|e^{\mathrm im\theta}|\leq e^{|m||\operatorname{Im}\theta|}$.
The polynomial norm, $f_N\leq32N$, and~\eqref{eq:uniform-majorant}
therefore imply
\[
 \|X_N\|+|p_N|
 \leq\frac{M}{N}
 \left(\frac{2\max\{|z|,|\bar z|\}
 e^{32|\operatorname{Im}\theta|}}{R_*}\right)^N.
\]
This gives normal convergence on every compact subset of the open set
\begin{equation}\label{eq:analytic-domain}
 2\max\{|z|,|\bar z|\}e^{32|\operatorname{Im}\theta|}<R_*,
 \qquad |\varepsilon|<2.
\end{equation}
Indeed, the ratio on any such compact subset is uniformly less than
one. Every term, including $X_0$, is holomorphic jointly in the
complex labels, $\theta$ and $\varepsilon$. The sums are therefore
jointly holomorphic on this open set, and all derivative series
converge uniformly on its compact subsets. For real $y$ we have
$|z|=|\bar z|=\|y\|$. Restriction to real variables proves joint
real analyticity on the common label neighbourhood
$\|y\|<R_*/2$, including $y=0$, for $\varepsilon\in[0,1]$ and
on a neighbourhood of both parameter endpoints. The dependence on
the real part of $\theta$ remains $2\pi$-periodic.

These convergence properties justify all products and derivatives in
the reduced equations. Every homogeneous coefficient of
$F_i-\partial_i p$ for $i=1,2$, and of $J-1$, vanishes by Proposition~\ref{prop:induction}. Their convergent sums consequently vanish identically,
so $F_i=\partial_i p$ and $J=1$ hold for the actual analytic functions
in the original local labels.

\section{Construction of the equilibria and their domains}\label{sec:construction}

In this section, we turn the convergent series of Section~\ref{sec:majorant}
into a family of analytic equilibria on solid tori.
Subsection~\ref{subsec:embedding} shows that $X$ embeds a common
neighbourhood of the axis and defines fields $B$ and $p$ that satisfy
\eqref{eq:mhs}. Subsection~\ref{subsec:unit-disk} selects the domain
$\Omega_\varepsilon$ as a pressure superlevel set by means of a Morse
chart and defines the unit-disk embedding $X_{\mathbb D}$. Together, these
give the solid tori and the analytic solutions of
Theorem~\ref{thm:main}, whose properties are derived in
Section~\ref{sec:properties}.

\subsection{Embedding and physical regularity}\label{subsec:embedding}

In this subsection, we show that $X$ embeds a common neighbourhood of
the axis for every $\varepsilon\in[0,1]$ and that its inverse gives
analytic fields $B$ and $p$ in physical coordinates, which satisfy
\eqref{eq:mhs}.

We use the original local labels and display the dependence on
$\varepsilon$. Section~\ref{sec:majorant} gives joint real analyticity
of $X(y,\theta;\varepsilon)$ and $p(y;\varepsilon)$ for
\[
 y\in U_0:=\{y\in\mathbb R^2:\|y\|<R_*/2\},
 \qquad \theta\in\mathbb S^1,\qquad \varepsilon\in(-2,2),
\]
together with $J=1$ and $F_i=\partial_i p$, $i=1,2$.
The joint map
\[
 \begin{aligned}
 \Psi:U_0\times\mathbb S^1\times(-2,2)
 &\longrightarrow\mathbb R^3\times(-2,2),\\
 (y,\theta,\varepsilon)
 &\longmapsto(X(y,\theta;\varepsilon),\varepsilon)
 \end{aligned}
\]
has invertible derivative because $J=1$. The real-analytic inverse
function theorem \cite[Section~2.5]{KrantzParks2002} therefore gives
local analytic inverses.

Set $I=[-1/2,3/2]$, so $[0,1]\subset\operatorname{int}I$.
We claim that there is an open disk $U\subset U_0$ centred at zero such
that $X(\,\cdot\,,\,\cdot\,;\varepsilon)$ is injective on
$U\times\mathbb S^1$ for every $\varepsilon\in I$.
Suppose, for contradiction, that no such disk exists. Then, for the corresponding sequences, there are
$\varepsilon_n\in I$ and distinct pairs
$(y_n,\theta_n),(\tilde y_n,\tilde\theta_n)\in U_0\times\mathbb S^1$
with $y_n,\tilde y_n\to0$ as \(n \to \infty\) and
\[
 X(y_n,\theta_n;\varepsilon_n)
 =X(\tilde y_n,\tilde\theta_n;\varepsilon_n).
\]
By compactness, pass to a subsequence with
$(\varepsilon_n,\theta_n,\tilde\theta_n)
\to(\varepsilon_*,\theta_*,\tilde\theta_*)
\in I\times\mathbb S^1\times\mathbb S^1$.
Since $X(0,\theta;\varepsilon)=e_1(\theta)$, taking limits gives
$e_1(\theta_*)=e_1(\tilde\theta_*)$. Injectivity of $e_1$ on
$\mathbb S^1$ gives $\theta_*=\tilde\theta_*$. The triples
$(y_n,\theta_n,\varepsilon_n)$ and $(\tilde y_n,\tilde\theta_n,\varepsilon_n)$
then both converge to $(0,\theta_*,\varepsilon_*)$. Since
$\varepsilon_*\in I\subset(-2,2)$, the inverse function theorem makes
$\Psi$ injective on a neighbourhood of this point. For large $n$ both
triples lie in this neighbourhood. They are distinct and have the same
image under $\Psi$, which is a contradiction. This proves the claim.

Thus $\Psi$ restricts to an analytic diffeomorphism from
$U\times\mathbb S^1\times\operatorname{int}I$ onto the open set
$V:=\Psi(U\times\mathbb S^1\times\operatorname{int}I)$.
For every $\varepsilon\in[0,1]$, the map
$X(\,\cdot\,,\,\cdot\,;\varepsilon)$ is therefore an analytic
embedding on $U\times\mathbb S^1$. With the pressure convention
of Section~\ref{sec:coordinates}, the relations
\[
 \begin{aligned}
 B(X(y,\theta;\varepsilon);\varepsilon)
   &=\partial_\theta X(y,\theta;\varepsilon),\\
 p(X(y,\theta;\varepsilon);\varepsilon)
   &=p(y;\varepsilon)
 \end{aligned}
\]
define physical fields jointly real analytic in
$(x,\varepsilon)\in V$, including on the axis.
The identities $J=1$ and $F_i=\partial_i p$, $i=1,2$, imply that our example is an MHD equilibrium \eqref{eq:mhs} by Lemma~\ref{lem:reduction}.

The disk $U$ is only a working neighbourhood of the axis. For
$\varepsilon>0$, the cubic coefficient $p_3$ is not radial by
Proposition~\ref{prop:cubic}, so the boundary of
$X(U\times\mathbb S^1;\varepsilon)$ need not be a pressure surface. In
Subsection~\ref{subsec:unit-disk}, we use a Morse chart to select
the domain $\Omega_\varepsilon$ inside $X(U\times\mathbb S^1;\varepsilon)$
as a pressure superlevel set.

\begin{figure}[ht!]
  \centering
  \includegraphics[width=\textwidth]{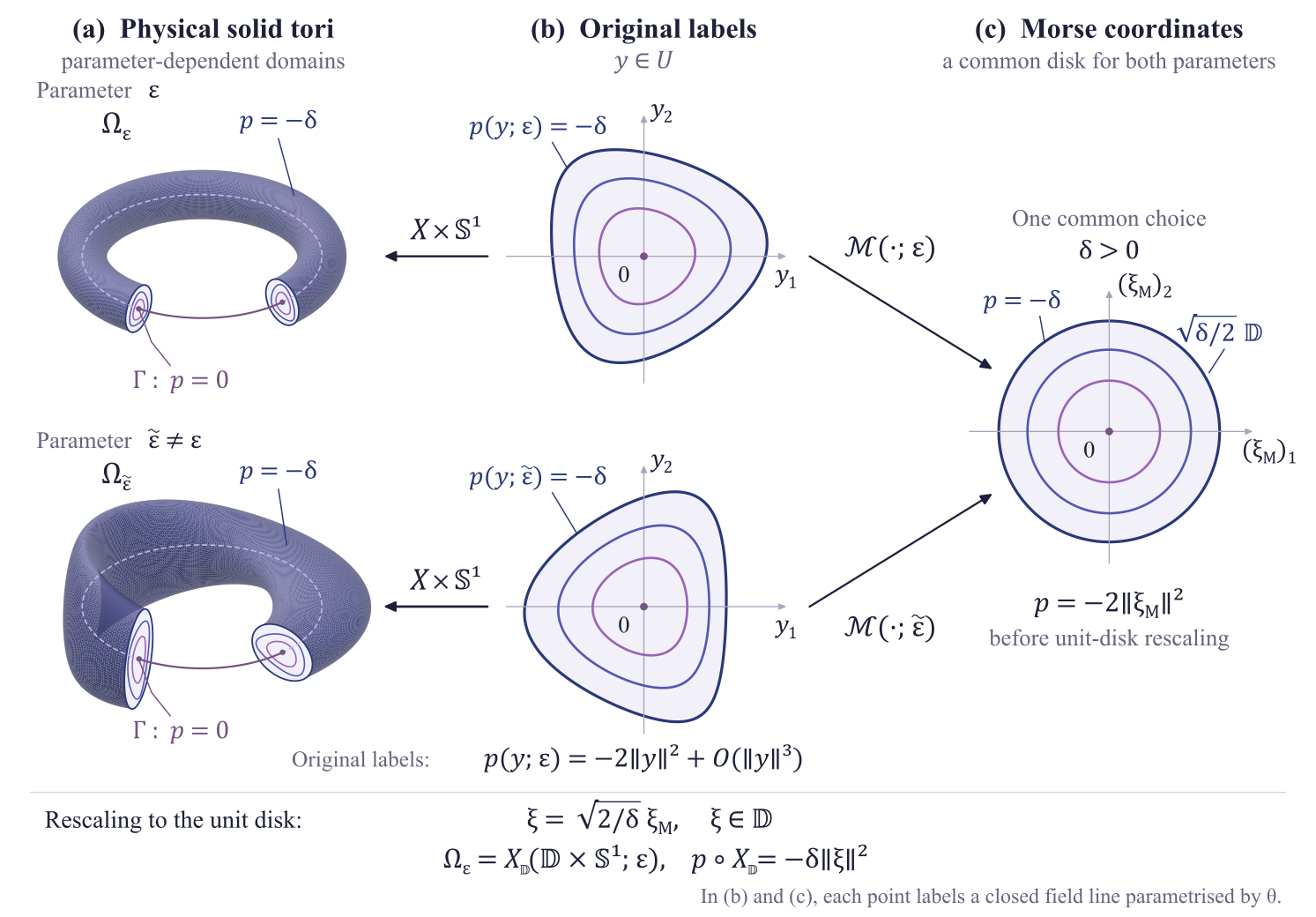}
  \caption{Schematic comparison of the domains for two distinct parameter
  values $\varepsilon,\tilde\varepsilon\in[0,1]$, using the same sufficiently
  small fixed $\delta>0$. First column: the corresponding physical solid
  tori, shown cut away to reveal the magnetic axis $\Gamma$ and nested
  pressure surfaces. Second column: the selected pressure superlevel
  components containing $y=0$ in the original label disk $U$.
  The products of these domains with $\mathbb S^1$ are mapped to the
    physical tori by the embedding $X$.
  Dark blue boundaries have $p=-\delta$; indigo and violet curves represent
  the same two selected pressure values between $-\delta$ and $0$ in both rows.
  Third column: the Morse map $\xi_{\mathrm M}=\mathcal M(y;\varepsilon)$
  sends the selected label domain to $\sqrt{\delta/2}\,\mathbb D$;
  the second row uses $\tilde\varepsilon$ in place of $\varepsilon$.
  In both cases the pressure equals $-2\|\xi_{\mathrm M}\|^2$, so the corresponding
  contours are the same concentric circles. The planar panels show only
  transverse coordinates: each point
  labels a closed field line parametrised by $\theta\in\mathbb S^1$,
  and the origin labels $\Gamma$. The final rescaling
  $\xi=\sqrt{2/\delta}\,\xi_{\mathrm M}$ gives
    $\Omega_\varepsilon=X_{\mathbb D}(\mathbb D\times\mathbb S^1;\varepsilon)$ with
  $p\circ X_{\mathbb D}=-\delta\|\xi\|^2$. The schematic shapes emphasise parameter
  dependence, they are not computed equilibria or drawn to scale.}
  \label{fig:domain-schematic}
\end{figure}

\subsection{Morse chart and the unit-disk embedding}\label{subsec:unit-disk}

In this subsection, we use the transverse pressure maximum at $y=0$
to construct a Morse chart, choose one boundary level $p=-\delta$ for
the whole family, and define the unit-disk embedding $X_{\mathbb D}$ and
the domains $\Omega_\varepsilon$.

By Proposition~\ref{prop:quadratic}, in the original local labels and
for every $\varepsilon\in I$,
\[
 p(y)=-2\|y\|^2+O(\|y\|^3),\qquad
 \partial_{i_1}\partial_{i_2}p(0)=-4\delta_{i_1i_2},
 \qquad i_1,i_2=1,2.
\]
Here and below $O(\cdot)$ denotes a Big-O, distinct from
the orthogonal matrix $O$ used for symmetries.

The pressure therefore has a nondegenerate transverse maximum at $y=0$.
We now construct transverse coordinates in which it is exactly quadratic.

Use the disk $U$ and parameter interval $I=[-1/2,3/2]$ from
Subsection~\ref{subsec:embedding}, shrinking $U$ when necessary.
In the original transverse labels, set
\[
 A(y;\varepsilon):=-\frac12\int_0^1(1-t)
       (D_y^2p)(ty;\varepsilon)\,dt,
 \qquad y\in U,\quad\varepsilon\in I.
\]
The matrix $A$ is symmetric, because $D_y^2p$ is. Since
$D_y^2p(0;\varepsilon)=-4\,\mathrm{Id}$
we have $A(0;\varepsilon)=\mathrm{Id}$ for every $\varepsilon\in I$.
Joint analyticity and compactness of $I$ allow $U$ to be chosen so
that $A$ is positive definite and its positive symmetric square
root is jointly real analytic. Define
\[
 \mathcal M(y;\varepsilon):=A(y;\varepsilon)^{1/2}y.
\]
Then $\mathcal M(0;\varepsilon)=0$ and
$D_y\mathcal M(0;\varepsilon)=\mathrm{Id}$. Since $p_0=p_1=0$, we have
$p(0;\varepsilon)=0$ and $\nabla_yp(0;\varepsilon)=0$. Taylor's formula
with integral remainder therefore gives
\[
 p(y;\varepsilon)
 =\int_0^1(1-t)\,y^{\mathsf T}(D_y^2p)(ty;\varepsilon)\,y\,dt
 =-2\,y^{\mathsf T}A(y;\varepsilon)\,y.
\]
Since $A^{1/2}$ is symmetric, $y^{\mathsf T}Ay=\|A^{1/2}y\|^2$, and hence
\begin{equation}\label{eq:morse-chart-norm}
 p(y;\varepsilon)=-2\|\mathcal M(y;\varepsilon)\|^2.
\end{equation}

By joint continuity and compactness of $I$, we may shrink $U$ so that
$\|D_y\mathcal M(y;\varepsilon)-\mathrm{Id}\|_{\mathrm{op}}<1/2$
for all $(y,\varepsilon)\in U\times I$, where $\|\cdot\|_{\mathrm{op}}$
denotes the operator norm induced by the Euclidean norm on $\mathbb R^2$.
Integrating along the segment joining $y,\tilde y\in U$ gives
\[
 \|\mathcal M(y;\varepsilon)-\mathcal M(\tilde y;\varepsilon)\|
 \geq\tfrac12\|y-\tilde y\|,
 \qquad\varepsilon\in I.
\]
Thus each $\mathcal M(\,\cdot\,;\varepsilon)$ is injective and the derivative bound also gives $\det D_y\mathcal M>0$.
The joint map
$(y,\varepsilon)\mapsto(\mathcal M(y;\varepsilon),\varepsilon)$
is therefore injective and has invertible derivative on
$U\times\operatorname{int}I$. Its image
\[
 \mathcal W:=\{(\mathcal M(y;\varepsilon),\varepsilon):
          y\in U,\ \varepsilon\in\operatorname{int}I\}
 \subset\mathbb R^2\times\operatorname{int}I
\]
is open, and its inverse
$(\xi_{\mathrm M},\varepsilon)\mapsto(\mathcal M^{-1}(\xi_{\mathrm M};\varepsilon),\varepsilon)$
is jointly real analytic on $\mathcal W$.
The coordinates $\xi_{\mathrm M}=\mathcal M(y;\varepsilon)$ form a Morse chart:
by \eqref{eq:morse-chart-norm}, the pressure in these coordinates
is exactly $-2\|\xi_{\mathrm M}\|^2$.
\par

The compact set $\{0\}\times[0,1]$ lies in $\mathcal W$. Choose one $\delta>0$, independent of $\varepsilon$, small enough for the winding estimates of Subsection~\ref{subsec:field-lines}, and such that
\[
 (\sqrt{\delta/2}\,\mathbb D)\times[0,1]\subset\mathcal W,
 \qquad \mathbb D=\{\xi\in\mathbb R^2:\|\xi\|\leq1\}.
\]

Define the unit-disk embedding by
\[
 \begin{aligned}
 X_{\mathbb D}(\xi,\theta;\varepsilon)
 &:=X\!\left(
      \mathcal M^{-1}(\sqrt{\delta/2}\,\xi;\varepsilon),
      \theta;\varepsilon\right),\\
 &(\xi,\theta,\varepsilon)\in
       \mathbb D\times\mathbb S^1\times[0,1].
 \end{aligned}
\]
Here $\xi=\sqrt{2/\delta}\,\xi_{\mathrm M}$ is the unit-disk label. Since the selected compact
set lies in the open inverse domain $\mathcal W$, this composition
is jointly real analytic on an open neighbourhood of
$\mathbb D\times\mathbb S^1\times[0,1]$.
The relabelling leaves $\theta$ unchanged, so
\[
 \begin{aligned}
 B(X_{\mathbb D}(\xi,\theta;\varepsilon);\varepsilon)
   &=\partial_\theta X_{\mathbb D}(\xi,\theta;\varepsilon),\\
 p(X_{\mathbb D}(\xi,\theta;\varepsilon);\varepsilon)
   &=-\delta\|\xi\|^2.
 \end{aligned}
\]

For each $\varepsilon\in[0,1]$, define
\[
 \Omega_\varepsilon=X_{\mathbb D}(\mathbb D\times\mathbb S^1;\varepsilon),
 \qquad
 \partial\Omega_\varepsilon
   =X_{\mathbb D}(\partial\mathbb D\times\mathbb S^1;\varepsilon).
\]
Since $X_{\mathbb D}$ is an analytic embedding, $\Omega_\varepsilon$ is a compact
solid torus with analytic boundary $p=-\delta$. It contains the common
axis
\[
 \Gamma=\{X_{\mathbb D}(0,\theta;\varepsilon):\theta\in\mathbb S^1\}
        =e_1(\mathbb S^1).
\]
Every original local label $y$ of $\Omega_\varepsilon$ has the form
$y=\mathcal M^{-1}(\sqrt{\delta/2}\,\xi;\varepsilon)$ with
$\xi\in\mathbb D$, so $\|\mathcal M(y;\varepsilon)\|\leq\sqrt{\delta/2}$.
The estimate
$\|\mathcal M(y;\varepsilon)-\mathcal M(\tilde y;\varepsilon)\|\geq\tfrac12\|y-\tilde y\|$
with $\tilde y=0$ therefore gives
\begin{equation}\label{eq:label-bound}
 \|y\|\leq\sqrt{2\delta} \, .
\end{equation}
We write $\Omega$ when $\varepsilon$ is fixed.

Figure~\ref{fig:domain-schematic} compares the physical and original
label domains for two parameter values. The original coordinates $y$
in panel~(b) are mapped to the Morse coordinates
$\xi_{\mathrm M}=\mathcal M(y;\varepsilon)$ in panel~(c). Rescaling this
common Morse disk to $\mathbb D$ gives the fixed unit-disk parametrisation.

\section{Properties of the solution}\label{sec:properties}

In this section, we derive the properties of the equilibria constructed
in Section~\ref{sec:construction}.
Subsections~\ref{subsec:pressure-surfaces}--\ref{subsec:asymmetry}
establish parts~\ref{item:part1}--\ref{item:part5} of Theorem~\ref{thm:main}, and
Subsection~\ref{sec:family} proves Theorem~\ref{thm:nonisolation}.

\subsection{Nonvanishing field and nested pressure surfaces}\label{subsec:pressure-surfaces}

In this subsection, we prove parts~\ref{item:part1} and~\ref{item:part2}
of Theorem~\ref{thm:main}.

Since $J=1$, the vector $\partial_\theta X$ does not vanish, so
$B\circ X=\partial_\theta X$ shows that $B$ is nowhere zero.
For $0<s\leq\delta$,
\[
 \{x\in\Omega_\varepsilon:p(x;\varepsilon)=-s\}
 =X_{\mathbb D}\!\left(
    \{\xi\in\mathbb D:\|\xi\|=\sqrt{s/\delta}\}
      \times\mathbb S^1;\varepsilon\right).
\]
These nested analytic tori foliate $\Omega_\varepsilon\setminus\Gamma$
and approach $\Gamma$ as $s\to0$. Since
$D_\xi(p\circ X_{\mathbb D})=-2\delta \xi$ and
$\partial_\theta(p\circ X_{\mathbb D})=0$, invertibility of $DX_{\mathbb D}$
shows that the critical set of physical pressure is exactly $\Gamma$.
Moreover,
\[
 (B\cdot\nabla p)\circ X_{\mathbb D}
 =\partial_\theta(p\circ X_{\mathbb D})=0,
\]
so $B$ is tangent to every pressure surface, including
$\partial\Omega_\varepsilon$. This proves parts~\ref{item:part1} and \ref{item:part2} of
Theorem~\ref{thm:main}. 

\subsection{Closed field lines and rotational transform}\label{subsec:field-lines}

In this subsection, we show that every field line is closed, that its
rotational transform away from $\Gamma$ is $2$, and that the domain
and equilibrium at $\varepsilon=0$ are axisymmetric.
Periodicity and the winding numbers give part~\ref{item:part3} of
Theorem~\ref{thm:main}, while the recursion at $\varepsilon=0$
yields the invariance under rotations about the vertical axis
asserted in part~\ref{item:part5}.

In the unit-disk labels, for each fixed $\xi\in\mathbb D$, the curve
$\theta\mapsto X_{\mathbb D}(\xi,\theta)$ is itself an integral curve of $B$,
because $B\circ X_{\mathbb D}=\partial_\theta X_{\mathbb D}$. Periodicity gives
period $2\pi$, and injectivity of $X_{\mathbb D}$ excludes every smaller
positive period.

Use the physical cylindrical coordinates
\[
 R=\sqrt{x_1^2+x_2^2},\qquad Z=x_3,\qquad
 x_1+\mathrm i x_2=Re^{\mathrm i\varphi}.
\]
Along an off-axis field line, the winding numbers about zero of
$(R-1)+\mathrm i Z$ and $x_1+\mathrm i x_2$ count its poloidal and
toroidal turns, respectively. Equivalently, one follows each
complex argument continuously over the field-line period without
reducing it modulo $2\pi$. This fixes our convention: the rotational
transform $\iota$ is the poloidal winding number divided by the
toroidal winding number.

In the original local labels $y=(y_1,y_2)$, with $z=y_1+iy_2$,
the prescribed linear term gives
\[
X=(1+a_1)e_1-a_2e_\theta+a_2e_2+O(\|y\|^2).
\]
Fix a closed disk $K\subset U_0$ centred at zero, independent of
$\delta$. Since $X$ is jointly
real analytic on a neighbourhood of the compact set
$K\times\mathbb S^1\times I$, Taylor's theorem in $y$ bounds the
remainder $X-X_0-X_1$ by a multiple of $\|y\|^2$ on $K$, with a constant
independent of $\theta$ and $\varepsilon$.
Hence $R-1=a_1+O(\|y\|^2)$ and $Z=a_2+O(\|y\|^2)$.
Since $a_1+ia_2=e^{2i\theta}z$, it follows that
\[
(R-1)+iZ=e^{2i\theta}z+O(\|y\|^2),
\qquad
x_1+ix_2=e^{i\theta}+O(\|y\|),
\]
uniformly in $\theta\in \mathbb{S}^1$ and $\varepsilon\in[0,1]$.
By \eqref{eq:label-bound}, every original local label $y$ of the
constructed torus satisfies $\|y\|\leq\sqrt{2\delta}$. By the choice of
$\delta$ in Subsection~\ref{subsec:unit-disk}, these labels lie in $K$,
and the estimates above give
\[
\bigl|(R-1)+iZ-e^{2i\theta}z\bigr|
<\frac{\|y\|}{2},
\qquad
\bigl|x_1+ix_2-e^{i\theta}\bigr|<\frac12
\]
for every nonzero original local label $y$ of the constructed torus
and every $\theta\in\mathbb S^1$.

For fixed $\varepsilon$ and $y\neq0$, consider the loops
\[
w_\lambda(\theta)
=(1-\lambda)e^{2i\theta}z
+\lambda\bigl((R-1)+iZ\bigr),
\qquad 0\leq\lambda\leq1,
\]
where $R$ and $Z$ are evaluated along the field line
$\theta\mapsto X(y,\theta)$.
The first bound implies
\[
|w_\lambda(\theta)|
\geq \|y\|
-\lambda\bigl|(R-1)+iZ-e^{2i\theta}z\bigr|
>\frac{\|y\|}{2}.
\]
Thus, normalisation by $|w_\lambda|$ gives a homotopy of maps
$w_\lambda/|w_\lambda|$ from $\mathbb{S}^1$ to $\mathbb{S}^1$.
By homotopy invariance of degree
\cite[Section~2.2, p.~134, property~(c)]{Hatcher2002},
the loop $(R-1)+iZ$ has the same winding number as
$e^{2i\theta}z$, namely $2$.
Applying the same argument to $x_1+ix_2$ and $e^{i\theta}$
gives toroidal winding number $1$.

Consequently, over one field-line period, the continuous increments
of the poloidal angle $\psi=\arg((R-1)+iZ)$ and the toroidal angle
$\varphi$ satisfy
\[
\Delta\psi=4\pi,
\qquad
\Delta\varphi=2\pi.
\]
The rotational transform, defined as the number of poloidal turns
per toroidal turn, is therefore
\[
\iota=\frac{\Delta\psi}{\Delta\varphi}=2
\]
on every off-axis field line, and hence on every regular pressure surface.



At $\varepsilon=0$, all free data vanish. Inductively, the frame
components of each $X_N$ are constant-coefficient polynomials in
the linear forms $a_1,a_2$, while each $p_N$ depends only on
$\|y\|$. This holds for the prescribed initial terms. To check the
induction step, differentiation at fixed $y$ uses
$\partial_\theta a_1=-2a_2$, $\partial_\theta a_2=2a_1$ and the
frame derivative identities, so it preserves this polynomial
form. The inverse $\Delta_y^{-1}$ commutes with label rotations,
because it preserves the angular weight of every
$z^j\bar z^k$. The transverse derivatives transform covariantly.
Angular averaging extracts a rotation-invariant polynomial, and
the zero-mean primitive $\partial_\theta^{-1}$ preserves the same
polynomial form for a mean-zero input. Thus
\eqref{eq:normal_preliminary}, \eqref{eq:normal_pressure},
\eqref{eq:pressure_recursion} and the tangential primitive in
\eqref{eq:embedding_recursion} preserve the induction claims.

Let $O$ be rotation through $\alpha$ about the vertical axis, and
let $\tilde y$ be the label with
$\tilde y_1+\mathrm i\tilde y_2=e^{-2\mathrm i\alpha}z$.
The linear forms $a_1,a_2$ are unchanged by simultaneously replacing
$(y,\theta)$ with $(\tilde y,\theta+\alpha)$, while the frame is
rotated by $O$. After convergence, the polynomial identities give
\begin{equation}\label{eq:epsilon-zero-equivariance}
 X(\tilde y,\theta+\alpha)=OX(y,\theta),\qquad p(\tilde y)=p(y).
\end{equation}
The latter identity preserves the local pressure superlevel component
used to construct $\Omega$. The Taylor matrix and its positive square
root transform covariantly under these label rotations at
$\varepsilon=0$, so the Morse chart and its inverse commute with them.
Thus \eqref{eq:epsilon-zero-equivariance} also holds with $X_{\mathbb D}$
and the unit-disk labels $\xi$ in place of $X$ and $y$. Differentiating the former identity at fixed labels and using $B\circ X=\partial_\theta X$ proves
$B(Ox)=OB(x)$, and the pressure convention proves $p(Ox)=p(x)$.
Thus the domain and equilibrium at $\varepsilon=0$ are invariant
under every rotation about the vertical axis.

\subsection{Asymmetry}\label{subsec:asymmetry}

We prove that every constructed equilibrium with $\varepsilon>0$
has no nonidentity Euclidean symmetry. First, the cubic pressure
coefficient determines third derivatives of physical pressure on
the axis. These exclude nonidentity isometries of $\Omega$ that
preserve $p$. Force balance then shows that a symmetry of $B$,
including one reversing its sign, would also preserve pressure.

Fix $\theta_0\in\mathbb S^1$ and
$v=v_1e_1(\theta_0)+v_2e_2$, with $v_1,v_2\in\mathbb R$.
This is a fixed normal vector at $e_1(\theta_0)$. Here $D$ denotes
differentiation in physical Cartesian coordinates, and $D^3p_x$
is the symmetric trilinear third derivative at $x$. Thus
\[
 D^3p_{e_1(\theta_0)}(v,v,v)
 =\left.\frac{d^3}{d\lambda^3}
             p(e_1(\theta_0)+\lambda v)\right|_{\lambda=0}.
\]
Both $\theta_0$ and $v$ are fixed in this straight-line derivative.
Define its inverse coordinates in the original local labels by
\[
 X(y(\lambda),\theta(\lambda))=e_1(\theta_0)+\lambda v,
 \qquad y(0)=0,\quad\theta(0)=\theta_0.
\]
Dots in this calculation denote derivatives with respect to
$\lambda$, so
\[
 y(\lambda)=\lambda\dot y(0)+\tfrac12\lambda^2\ddot y(0)
                      +O(\lambda^3).
\]
Differentiating the inverse-coordinate identity once gives
\begin{equation}\label{eq:asymmetry-inverse-first}
 \dot y(0)=
 (v_1\cos2\theta_0+v_2\sin2\theta_0,
       -v_1\sin2\theta_0+v_2\cos2\theta_0),
 \qquad \dot\theta(0)=v_2.
\end{equation}
In the next formulas, $X_1,X_2$ and their angular derivatives
are evaluated at $\theta_0$; their first arguments remain labels.
The second derivative of that identity is
\begin{equation}\label{eq:asymmetry-inverse-second}
\begin{aligned}
 0={}&X_1(\ddot y(0),\theta_0)+2X_2(\dot y(0),\theta_0)
       +2v_2\partial_\theta X_1(\dot y(0),\theta_0)\\
    &-v_2^2e_1(\theta_0)+\ddot\theta(0)e_\theta(\theta_0).
\end{aligned}
\end{equation}
Here
$\partial_\theta X_1(\dot y(0),\theta_0)
 =-v_2e_1(\theta_0)-v_1e_\theta(\theta_0)+2v_1e_2$,
with $\partial_\theta$ still taken at fixed labels. Also
$v\cdot X_1(\ddot y(0),\theta_0)=\dot y(0)\cdot\ddot y(0)$
and $v\cdot e_\theta(\theta_0)=0$. Taking the scalar product of
\eqref{eq:asymmetry-inverse-second} with $v$ therefore gives
\begin{equation}\label{eq:asymmetry-inverse-contraction}
 \dot y(0)\cdot\ddot y(0)
 =-2v\cdot X_2(\dot y(0),\theta_0)-v_1v_2^2.
\end{equation}
The last term includes the variation of the angular coordinate
along the physical straight line.

Since $p(y)=-2\|y\|^2+p_3(y)+O(\|y\|^4)$, the chain rule gives
\begin{equation}\label{eq:asymmetry-cubic-chain-rule}
\begin{aligned}
 D^3p_{e_1(\theta_0)}(v,v,v)
 &=6p_3(\dot y(0))-12\dot y(0)\cdot\ddot y(0)\\
 &=6p_3(\dot y(0))+24v\cdot X_2(\dot y(0),\theta_0)+12v_1v_2^2.
\end{aligned}
\end{equation}
At these labels, $a_1=v_1$, $a_2=v_2$ and
$\|\dot y(0)\|^2=v_1^2+v_2^2$. By Proposition~\ref{prop:quadratic},
\[
 v\cdot X_2(\dot y(0),\theta_0)
 =-\tfrac12v_1(v_1^2+v_2^2)
   +\frac{\varepsilon g(\theta_0)}{2}(v_1^3-3v_1v_2^2).
\]
Substitution of this identity and Proposition~\ref{prop:cubic} yields
\begin{equation}\label{eq:asymmetry-cubic-normal}
\begin{aligned}
 D^3p_{e_1(\theta_0)}(v,v,v)
 = &-12v_1^3+12\varepsilon g(\theta_0)(v_1^3-3v_1v_2^2)
 +7\varepsilon^2
 (v_2\cos2\theta_0-v_1\sin2\theta_0)(v_1^2+v_2^2).
\end{aligned}
\end{equation}
These derivatives belong to the full analytic pressure. Terms
$p_N$ with $N\geq4$ and $X_N$ with $N\geq3$ cannot enter:
$\nabla_y p(0)=0$, so the cubic chain rule uses only $p_2,p_3$
and the first two derivatives of the inverse map. Convergence of
the derivative series justifies their evaluation on the axis.
In particular, as functions of the axis angle,
\begin{align}
 D^3p_{e_1(\theta)}(e_1(\theta),e_1(\theta),e_1(\theta))
 &=-12+12\varepsilon(\cos3\theta+\sin5\theta)
       -7\varepsilon^2\sin2\theta,
 \label{eq:asymmetry-radial-derivative}\\
 D^3p_{e_1(\theta)}(e_2,e_2,e_2)
 &=7\varepsilon^2\cos2\theta.
 \label{eq:asymmetry-vertical-derivative}
\end{align}

Suppose $G(x)=Ox+\tau$, $O\in\mathrm O(3)$, maps $\Omega$ onto
itself and $p\circ G=p$. It preserves the critical set $\Gamma$.
An isometry preserving this circle preserves its centre, its
affine plane and its normal line. Therefore
\begin{equation}\label{eq:asymmetry-circle-isometry}
 \tau=0,\qquad O=O_{\mathrm h}\oplus\sigma,
 \qquad O_{\mathrm h}\in\mathrm O(2),\quad \sigma\in\{1,-1\}.
\end{equation}
The horizontal action on the circle is either
$\theta\mapsto\theta+\alpha$ or
$\theta\mapsto-\theta+\alpha$. It sends each radial vector to
the radial vector at the image point. Differentiating $p\circ G=p$
three times thus makes the radial derivative in
\eqref{eq:asymmetry-radial-derivative} invariant under this action.

For a rotation, its nonzero third and fifth Fourier components require
\begin{equation}\label{eq:asymmetry-rotation}
 3\alpha=0\pmod{2\pi},\qquad5\alpha=0\pmod{2\pi}.
\end{equation}
Since $3$ and $5$ are coprime, $\alpha=0\pmod{2\pi}$.
For a reflection the same components require instead
\begin{equation}\label{eq:asymmetry-reflection}
 3\alpha=0\pmod{2\pi},\qquad5\alpha=\pi\pmod{2\pi}.
\end{equation}
Writing $3\alpha=2\pi j$ and $5\alpha=(2k+1)\pi$ for integers
$j,k$ would give $10j=3(2k+1)$, an even integer equal to an odd
integer. Thus reflection is impossible. Distinct Fourier
frequencies cannot cancel, so every nonidentity rotation and
every reflection of the circle is excluded. Hence
$O_{\mathrm h}=\mathrm{Id}$.

The only remaining nonidentity possibility is
$G_S(x)=(x_1,x_2,-x_3)$. It fixes $e_1(0)$ and sends $e_2$ to
$-e_2$. Pressure invariance would make
$D^3p_{e_1(0)}(e_2,e_2,e_2)$ equal to its negative, whereas
\eqref{eq:asymmetry-vertical-derivative} gives its nonzero value
$7\varepsilon^2$. Thus $p$ has no nonidentity Euclidean
symmetry of $\Omega$.

Finally suppose $G(x)=Ox+\tau$ maps $\Omega$ onto itself and
\[
 B\circ G=\eta OB,\qquad O\in\mathrm O(3),\qquad\eta\in\{1,-1\}.
\]
The curl transformation law is
\[
 (\nabla\times B)\circ G=\eta\det(O)\,O(\nabla\times B).
\]
Using the cross-product transformation under an orthogonal matrix,
\[
 [ (\nabla\times B)\times B]\circ G
 =\eta^2(\det O)^2O[(\nabla\times B)\times B]
 =O[(\nabla\times B)\times B].
\]
The chain rule and force balance at $Gx$ and $x$ therefore give
\[
 \nabla(p\circ G)=O^{\mathsf T}(\nabla p)\circ G,
 \qquad (\nabla p)\circ G=O\nabla p,
\]
and hence $\nabla(p\circ G)=\nabla p$.
Connectedness of $\Omega$ implies $p\circ G=p+c$ for a constant
$c$. Since $\Omega$ is compact and $G$ maps it onto itself,
$\max_\Omega(p\circ G)=\max_\Omega p$, forcing $c=0$.
The pressure argument now gives $G=\mathrm{Id}$, and the
nonvanishing of $B$ gives $\eta=1$. This completes the proof of Theorem \ref{thm:main}.

\subsection{Nonisolation and the role of the boundary}\label{sec:family}

\begin{proof}[Proof of Theorem~\ref{thm:nonisolation}]
By Subsection~\ref{subsec:unit-disk}, the embedding $X_{\mathbb D}$ has
the same domain $\mathbb D\times\mathbb S^1$ for every
$\varepsilon\in[0,1]$ and is jointly real analytic on an open
neighbourhood of
$\mathbb D\times\mathbb S^1\times[0,1]$. In these unit-disk labels,
\[
 B\circ X_{\mathbb D}=\partial_\theta X_{\mathbb D} \, ,\qquad
 p\circ X_{\mathbb D}=-\delta\|\xi\|^2 \, .
\]
The first composition is jointly real analytic because $X_{\mathbb D}$ is,
and the second is independent of $\varepsilon$ because $\delta$
was chosen once for the whole family. This proves the analyticity
assertions, including at both parameter endpoints.

Distinct members of the family are not related by operations that
trivially preserve \eqref{eq:mhs}. By
\eqref{eq:asymmetry-radial-derivative}, the function
$\theta\mapsto D^3p_{e_1(\theta)}(e_1(\theta),e_1(\theta),e_1(\theta))$
has mean $-12$, and its Fourier coefficients of frequency $\pm3$ have
modulus $6\varepsilon$. Their ratio $\varepsilon/2$, computed on the
magnetic axis, is unchanged by Euclidean similarities, by multiplying
$B$ by a nonzero constant and $p$ by its square, by adding a constant
to $p$, by reparametrising $\mathbb D\times\mathbb S^1$, and by
restriction to a smaller domain bounded by a pressure surface, since it
depends only on the third derivatives of $p$ on $\Gamma$. Combined with the joint analyticity just proved, every member is
therefore the $C^\infty$ limit of members not related to it by these
operations, which proves the non-isolation assertion. In particular, the
axisymmetric equilibrium at $\varepsilon=0$ is such a limit of equilibria
with trivial symmetry group.


The family $\varepsilon \mapsto (X_{\mathbb D},B\circ X_{\mathbb D},p\circ X_{\mathbb D})$ therefore provides nonisolated (with respect to the \(C^\infty \) topology) equilibrium--domain
pairs and a branch attached to an axisymmetric equilibrium.
Its boundary is the pressure surface $p=-\delta$ for each
parameter, parametrised by the fixed $\partial\mathbb D\times\mathbb S^1$.
\end{proof}

\providecommand{\MHSVisOrder}{18}
\providecommand{\MHSVisOuterResidual}{3.62\%}
\providecommand{\MHSVisFigureDirectory}{figures}

\begin{figure}[th]
 \centering
 \includegraphics[width=\textwidth]{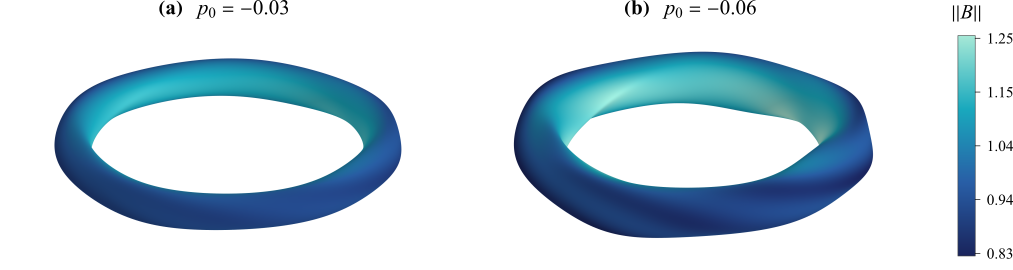}
 \caption{Pressure surfaces of the degree-$\MHSVisOrder$ truncation
 with $\varepsilon=1$:
 (a) $p_0=-0.03$ and (b) $p_0=-0.06$. Colour gives $\|B\|$ with the
 same scale in both panels. Spatial aspect ratios are preserved.}
 \label{fig:mhs-pressure-surfaces}
\end{figure}

\subsection{Visualisation}\label{sec:figures}

\begin{figure}[th]
 \centering
 \includegraphics[width=\textwidth]{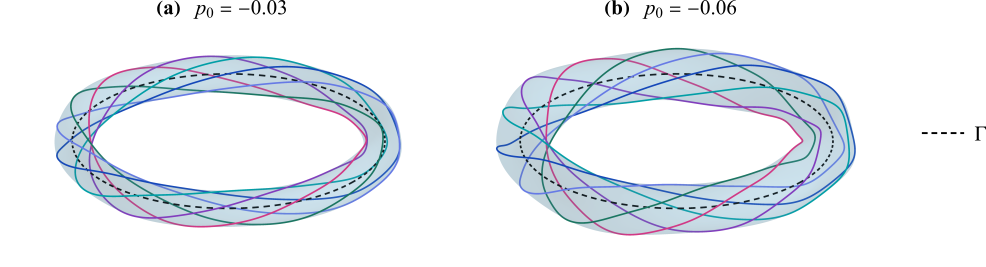}
 \caption{Six selected field lines on each pressure surface in
 Figure~\ref{fig:mhs-pressure-surfaces}: (a) $p_0=-0.03$ and
 (b) $p_0=-0.06$. Each coloured curve is the closed orbit
 $\theta\mapsto X_{<N+1}(y,\theta)$ at fixed transverse labels $y$.
 The colours identify the same orbits in
 Figure~\ref{fig:mhs-poincare}. The dashed curve is the magnetic axis
 $\Gamma$.}
 \label{fig:mhs-field-lines}
\end{figure}
\begin{figure}[th]
 \centering
 \includegraphics[width=\textwidth]{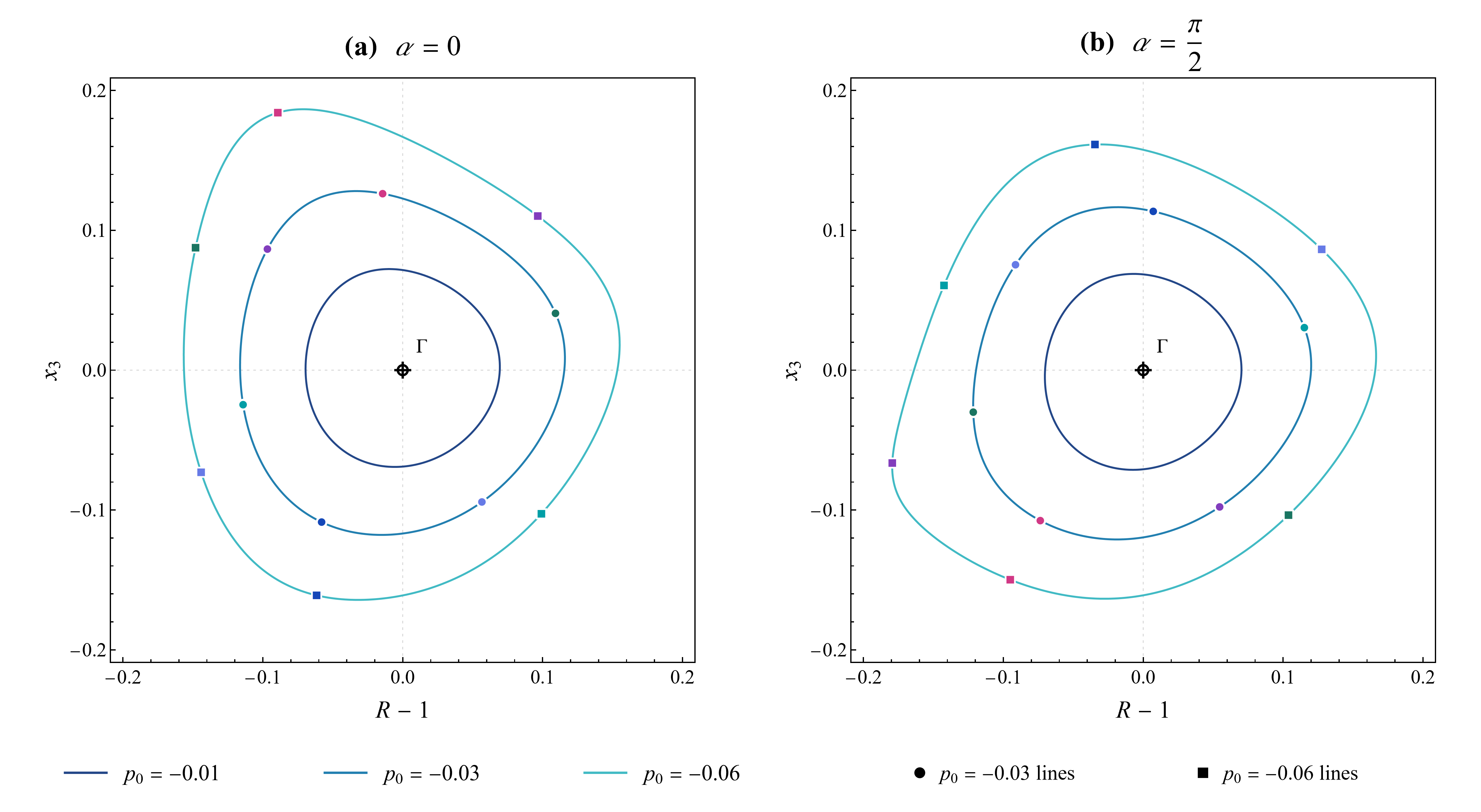}
 \caption{Physical normal sections at (a) $\alpha=0$ and
 (b) $\alpha=\pi/2$. Curves show intersections of the pressure
 surfaces labelled $p_0=-0.01,-0.03,-0.06$ with the section in
 \eqref{eq:visualization-sections}.
 Coloured punctures correspond to the field lines in
 Figure~\ref{fig:mhs-field-lines}; circles mark $p_0=-0.03$ and
 squares mark $p_0=-0.06$. Since every orbit is closed, repeated
 transits give the same puncture.
 The marked origin is the magnetic-axis intersection with the section.
 Both axes are physical coordinates, $(R-1,x_3)$, with equal
 length scales.}
 \label{fig:mhs-poincare}
\end{figure}

Figures~\ref{fig:mhs-pressure-surfaces}--\ref{fig:mhs-poincare}
illustrate finite truncations of the construction. The figures use
$\varepsilon=1$. No convergence-radius bound has been certified
for the plotted pressure levels.

For a truncation degree $N\geq 2$, use the finite sums
\begin{equation}\label{eq:visualization-truncations}
 X_{<N+1}=\sum_{j=0}^N X_j,\qquad
 p_{<N+1}(y)=\sum_{j=2}^N p_j(y).
\end{equation}
Where the truncated map is invertible, define its physical fields by
\begin{equation}\label{eq:visualization-fields}
 B_{<N+1}\circ X_{<N+1}=\partial_\theta X_{<N+1},\qquad
 p_{<N+1}(X_{<N+1}(y,\theta))=p_{<N+1}(y),
\end{equation}
using the same physical/label pressure convention as in \eqref{eq:BP}.
The figures use $N=\MHSVisOrder$, while the residual comparison below
also includes lower truncation degrees.
Throughout this subsection, $p$ denotes the pressure field and
$p_0<0$ denotes a prescribed pressure value. Each plotted
surface satisfies $p_{<N+1}(x)=p_0$.
The displayed levels are $p_0=-0.01,-0.03,-0.06$.
Truncation subscripts are suppressed in the figure captions.
All spatial views preserve Euclidean aspect ratios.

Figure~\ref{fig:mhs-pressure-surfaces} shows the two larger surfaces,
with $p_0=-0.03$ and $p_0=-0.06$, coloured by magnetic-field strength
on one common scale. Figure~\ref{fig:mhs-field-lines} shows six
selected field lines on each surface. Each curve holds $y$ fixed
and varies $\theta$ over one period, so its tangent is
$B_{<N+1}$ by \eqref{eq:visualization-fields}.
The axis $\Gamma=e_1(\mathbb S^1)$ is included for reference.

The sections in Figure~\ref{fig:mhs-poincare} are physical normal
sections of the circular axis. For an angle $\alpha$, the section is
\begin{equation}\label{eq:visualization-sections}
 \{x\in\mathbb R^3:
  (x-e_1(\alpha))\cdot e_\theta(\alpha)=0,\quad
  x\cdot e_1(\alpha)>0\},
\end{equation}
restricted to the local component near $e_1(\alpha)$.
The displayed coordinates are $(R-1,x_3)$, where
$R=x\cdot e_1(\alpha)$ and $x_3=x\cdot e_2$ on this section.
The angles are $\alpha=0,\pi/2$. Intersections are computed with
the physical plane \eqref{eq:visualization-sections}, rather than
by setting the field-line coordinate $\theta$ equal to $\alpha$.
The three pressure contours and intersections of the field lines
in Figure~\ref{fig:mhs-field-lines} are shown. Every field line
closes after one $2\pi$ period in $\theta$, with rotational
transform $2$. Repeated returns therefore coincide; only distinct
punctures are plotted. The marked origin is the intersection of
$\Gamma$ with the section.

The finite-order force balance is assessed through the physical
relative residual
\begin{equation}\label{eq:visualization-residual}
 \mathcal E_N(x)=
 \frac{\bigl\|(\nabla\times B_{<N+1})\times B_{<N+1}
                         -\nabla p_{<N+1}\bigr\|}
      {\|\nabla p_{<N+1}\|}.
\end{equation}
All derivatives and norms in this expression are taken in physical
Euclidean coordinates. We compare truncation degrees
$N=2,4,\ldots,18$ at the three prescribed pressure values
$p_0=-0.01,-0.03,-0.06$, keeping $\varepsilon=1$ fixed.
Here $N$ denotes the largest homogeneous degree retained in both
the map and the pressure.

For each pair $(N,p_0)$, the pressure surface is recomputed from
the corresponding truncated pressure. Writing
\[
 y=r_N(\beta;p_0)(\cos\beta,\sin\beta),\qquad
 p_{<N+1}(y)=p_0,
\]
we sample its image under $X_{<N+1}(y,\theta)$ on a uniform periodic
grid in $(\beta,\theta)$. If $\mathcal G_N(p_0)$ denotes this set
of sampled physical points, the reported quantity is
\begin{equation}\label{eq:visualization-sampled-residual}
 \widehat{\mathcal E}_N(p_0)
 =\max_{x\in\mathcal G_N(p_0)}\mathcal E_N(x).
\end{equation}
Table~\ref{tab:visualization-residual} gives
$100\widehat{\mathcal E}_N(p_0)$, expressed as a percentage, using
a $192\times1024$ grid. Comparison with the nested
$96\times512$ grid changes each sampled maximum by at most
$0.233\%$ of its value on the finer grid.

The pressure level provides a measure of proximity to the magnetic
axis. Near $\Gamma$,
\[
 p_{<N+1}(y)=-2\|y\|^2+O(\|y\|^3),\qquad
 \operatorname{dist}\!\left(X_{<N+1}(y,\theta),\Gamma\right)
   =\|y\|+O(\|y\|^2).
\]
Thus $\sqrt{-p_0/2}$ is the leading-order distance scale, although
a pressure surface is not generally a surface of constant Euclidean
distance from the axis.

\begin{table}[tbp]
 \centering
 \caption{Maximum sampled relative force-balance residual
 $100\widehat{\mathcal E}_N(p_0)$, in percent, for
 $\varepsilon=1$. The degree $N$ is the largest homogeneous
 degree retained in both the map and the pressure.
 Each pressure surface is recomputed at its own truncation degree
 and sampled on a $192\times1024$ periodic grid in
 $(\beta,\theta)$. These are sampled maxima, not certified
 upper bounds over the surfaces. Taking this data and fitting the log-error $ \log (100 \, \widehat{\mathcal E}_N(-0.06) )$ we obtain a fit \(5.69 - 0.238 N\) with \(r^2\) value \(0.995\), indicating that the corresponding error is decaying exponentially in the truncation degree.
 }
 \label{tab:visualization-residual}
 \begin{tabular}{rccc}
  \hline
  $N$ & $p_0=-0.01$ & $p_0=-0.03$ & $p_0=-0.06$ \\
  \hline
   2 & $56.8$                 & $115$                  & $187$  \\
   4 & $5.39$                 & $33.0$                 & $99.7$ \\
   6 & $0.607$                & $11.0$                 & $65.5$ \\
   8 & $7.95\times10^{-2}$    & $4.24$                 & $48.0$ \\
  10 & $9.21\times10^{-3}$    & $1.43$                 & $30.8$ \\
  12 & $1.03\times10^{-3}$    & $0.461$                & $18.8$ \\
  14 & $1.11\times10^{-4}$    & $0.144$                & $11.0$ \\
  16 & $1.23\times10^{-5}$    & $4.41\times10^{-2}$     & $6.35$ \\
  18 & $1.35\times10^{-6}$    & $1.34\times10^{-2}$     & $3.62$ \\
  \hline
 \end{tabular}
\end{table}

For the degrees and pressure levels tested, the sampled residual
increases with $-p_0$ and decreases roughly geometrically in $N$:
from $N=14$ to $N=18$, each increase of the degree by two reduces it
by a factor of about $9$, $3.3$ and $1.7$ at $p_0=-0.01$, $-0.03$
and $-0.06$.
At $N=18$, the residuals are approximately
$1.35\times10^{-6}\%$, $0.0134\%$, and $3.62\%$
at $p_0=-0.01,-0.03,-0.06$, respectively.
These diagnostics measure the remaining force imbalance of the
finite truncations. The sampled maxima and their grid-refinement
comparison provide neither a certified bound over an entire surface
nor a bound on the difference from the infinite series.
The figures and table illustrate the finite-order construction
and do not enter the proof.

\clearpage
\appendix

\section{Explicit coefficients of degrees two and three}\label{app:explicit}

In this appendix, we evaluate the recursion of Definition~\ref{def:example}
at degrees two and three. Subsection~\ref{subsec:unit-disk} uses
$p_2$, and Subsection~\ref{subsec:asymmetry} uses $X_2$, $p_2$ and $p_3$ to
compute the third derivatives of physical pressure along normal directions
at $\Gamma$ and thereby prove part~\ref{item:part4} of Theorem~\ref{thm:main}.

\begin{prop}[Quadratic coefficients]\label{prop:quadratic}
The coefficients $X_2$ and $p_2$ of Definition~\ref{def:example} are
given by
\begin{equation}\label{eq:deg2}
\begin{aligned}
 p_2&=-2\|y\|^2,\qquad
 X_2=-\tfrac12\|y\|^2e_1+w^{\mathrm f},\\
 w^{\mathrm f}&=w_1^{\mathrm f}e_1+w_\theta^{\mathrm f}e_\theta
                 +w_2^{\mathrm f}e_2,\\
 w_1^{\mathrm f}&=\tfrac{\varepsilon g}{2}(a_1^2-a_2^2),\qquad
 w_2^{\mathrm f}=-\varepsilon g\,a_1a_2,\\
 \partial_\theta w_\theta^{\mathrm f}&=-2w_1^{\mathrm f},\qquad
 \langle w_\theta^{\mathrm f}\rangle=0.
\end{aligned}
\end{equation}
\end{prop}
\begin{proof}
For $N=2$, the degree-one components of \eqref{eq:base_functionals} and \eqref{PiNi} give
\[
 \hat F_i=-6y_i,\qquad \mathcal{F}=-3\|y\|^2,\qquad
 \hat J=a_1,\qquad \hat W=0.
\]
It follows that
$\hat C=\partial_\theta^{-1}(4a_1-6a_1)=-a_2$.
Rotation invariance of the right inverse gives
\[
 \Delta_y^{-1}a_1=\tfrac18a_1\|y\|^2,\qquad
 \Delta_y^{-1}(-a_2)=-\tfrac18a_2\|y\|^2,
 \qquad
 \hat w=\left(-\tfrac14(a_1^2-a_2^2),-\tfrac12a_1a_2\right)+H_2.
\]
The first component has zero mean: its nonfree part has frequency four,
and its free part is a product of frequencies four with three or five.
Thus $p_2=\tfrac23\langle-3\|y\|^2-2\hat w_1\rangle=-2\|y\|^2$.
The potential and the planar pressure correction are
\[
 \Delta_y^{-1}p_2=-\tfrac18\|y\|^4,\qquad
 -\tfrac12\nabla_a^\perp\partial_{a_2}\Delta_y^{-1}p_2
 =\left(-\tfrac14(a_1^2+3a_2^2),\tfrac12a_1a_2\right).
\]
Adding them to $\hat w$ yields $X_2$ and $p_2$ as in \eqref{eq:deg2}.
The last line of \eqref{eq:deg2} follows from
$\mathcal{F}-p_2-2w_1=-2w_1^{\mathrm f}$; its right-hand side has zero
mean, so it uniquely specifies the stated periodic primitive.
\end{proof}

\begin{prop}[Cubic pressure coefficient]\label{prop:cubic}
The coefficient $p_3$ of Definition~\ref{def:example} is
\[
 p_3=\frac{7\varepsilon^2}{6}y_2\|y\|^2.
\]
\end{prop}
\begin{proof}
Throughout this proof $N=3$, and
$\hat F_i,\mathcal{F},\hat J,\hat W,\hat C,\hat w$ are the sources
and preliminary normal correction already defined for this degree.
By \eqref{eq:X1} and Proposition~\ref{prop:quadratic},
$X_1=a_1e_1-a_2e_\theta+a_2e_2$ and
$X_2=-\tfrac12\|y\|^2e_1+w^{\mathrm f}$, with $w^{\mathrm f}$ as in
\eqref{eq:deg2}.
Since $H_3=0$, the pressure recursion becomes
\begin{equation}\label{eq:Pi3mean}
 p_3=\frac23\left\langle\mathcal{F}-2\hat w_1\right\rangle,
 \qquad
 \hat w_1=-\partial_{a_1}\Delta_y^{-1}\hat J
          -\partial_{a_2}\Delta_y^{-1}\hat C.
\end{equation}
Let $\mathcal P_2$ retain only the Fourier frequencies $2$ and $-2$.
The coefficients in $\partial_{a_1}$ and $\partial_{a_2}$ have
these frequencies. Therefore the average in \eqref{eq:Pi3mean}
depends only on $\langle\mathcal{F}\rangle$,
$\mathcal P_2\hat J$ and $\mathcal P_2\hat C$.

Only the pairs $(X_1,X_2)$ and $(X_2,X_1)$ contribute to the
degree-two part of $F_i$. Homogeneity in $y$ gives
\[
 \mathcal{F}=-\partial_\theta X_1\cdot\partial_\theta X_2
       +\tfrac13\partial_\theta^2X_2\cdot X_1
       +\tfrac23\partial_\theta^2X_1\cdot X_2.
\]
Using the derivatives of $a_1,a_2$ and of the frame,
\[
 \partial_\theta X_1=-a_2e_1-a_1e_\theta+2a_1e_2,\qquad
 \partial_\theta^2X_1=-a_1e_1+a_2e_\theta-4a_2e_2.
\]
The contribution of $-\|y\|^2e_1/2$ to $\mathcal{F}$ is
\[
 (-\tfrac12+\tfrac16+\tfrac13)a_1\|y\|^2=0.
\]
The frame components of $X_1$ have frequency $2$, whereas
$w_1^{\mathrm f},w_2^{\mathrm f}$ and their zero-mean tangential
primitive have only odd frequencies. Differentiation in the
rotating frame preserves this parity. The remaining terms in
$\mathcal{F}$ therefore have only odd frequencies, giving
\[
 \langle\mathcal{F}\rangle=0,\qquad
 \mathcal P_2\partial_{a_1}\mathcal{F}=0.
\]

Expanding the degree-two part of the Jacobian gives
\[
\begin{aligned}
 \hat J={}&\det(\partial_1X_1,\partial_\theta X_2,\partial_2X_1)
          +\det(\partial_1X_2,\partial_\theta X_1,\partial_2X_1)\\
 &+\det(\partial_1X_1,\partial_\theta X_1,\partial_2X_2)
          +\det(\partial_1X_2,e_\theta,\partial_2X_2).
\end{aligned}
\]
Its contribution independent of $H_2$ is
$-\|y\|^2/2-a_1^2$, which has only frequencies $0,\pm4$.
The terms linear in $w^{\mathrm f}$ have only odd frequencies.
The quadratic contribution is
\[
 \det(\partial_1w^{\mathrm f},e_\theta,\partial_2w^{\mathrm f})
 =\det\frac{\partial(w_1^{\mathrm f},w_2^{\mathrm f})}
                 {\partial(y_1,y_2)}
 =-\varepsilon^2g^2\|y\|^2.
\]
Indeed, tangential components do not enter this determinant, the
linear forms $(a_1,a_2)$ are obtained from $y$ by a determinant-one
rotation, and their derivative matrix for the two normal components is
\[
 \varepsilon g
 \begin{pmatrix}a_1&-a_2\\-a_2&-a_1\end{pmatrix}.
\]
Since
\[
 g^2=1+\tfrac12\cos6\theta-\tfrac12\cos10\theta
                 +\sin8\theta+\sin2\theta,
\]
we obtain
\[
 \mathcal P_2\hat J=-\varepsilon^2\|y\|^2\sin2\theta.
\]

Only the pair $(X_2,X_2)$ contributes to the degree-two part of $W$:
\[
 \hat W=\partial_1X_2\cdot\partial_2\partial_\theta X_2
       -\partial_1\partial_\theta X_2\cdot\partial_2X_2.
\]
The contribution independent of $H_2$ vanishes and the linear
terms again have only odd frequencies. In the quadratic part,
the tangential and rotating-frame terms cancel:
\[
\begin{aligned}
 &(\partial_1w_\theta^{\mathrm f})(\partial_2\partial_\theta w_\theta^{\mathrm f})
 -(\partial_1\partial_\theta w_\theta^{\mathrm f})(\partial_2w_\theta^{\mathrm f})\\
 &\quad+2\bigl[(\partial_1w_\theta^{\mathrm f})(\partial_2w_1^{\mathrm f})
              -(\partial_1w_1^{\mathrm f})(\partial_2w_\theta^{\mathrm f})\bigr]=0,
\end{aligned}
\]
by $\partial_\theta w_\theta^{\mathrm f}=-2w_1^{\mathrm f}$.
For the remaining normal terms, write directly
\[
 w_1^{\mathrm f}+\mathrm i w_2^{\mathrm f}
 =\frac{\varepsilon g}{2}(a_1-\mathrm i a_2)^2
 =\frac{\varepsilon g}{2}e^{-4\mathrm i\theta}\bar z^2.
\]
Then
\[
\begin{aligned}
 \partial_1(w_1^{\mathrm f}+\mathrm i w_2^{\mathrm f})
   &=\varepsilon g\,e^{-4\mathrm i\theta}\bar z,\\
 \partial_2(w_1^{\mathrm f}+\mathrm i w_2^{\mathrm f})
   &=-\mathrm i\,\partial_1(w_1^{\mathrm f}+\mathrm i w_2^{\mathrm f}),\\
 \overline{\partial_1(w_1^{\mathrm f}+\mathrm i w_2^{\mathrm f})}\,
 \partial_\theta\partial_1(w_1^{\mathrm f}+\mathrm i w_2^{\mathrm f})
   &=\varepsilon^2\|y\|^2
                  (g\partial_\theta g-4\mathrm i g^2).
\end{aligned}
\]
Consequently,
\[
\begin{aligned}
 &(\partial_1w_1^{\mathrm f})(\partial_2\partial_\theta w_1^{\mathrm f})
  -(\partial_1\partial_\theta w_1^{\mathrm f})(\partial_2w_1^{\mathrm f})\\
 &\quad+(\partial_1w_2^{\mathrm f})(\partial_2\partial_\theta w_2^{\mathrm f})
  -(\partial_1\partial_\theta w_2^{\mathrm f})(\partial_2w_2^{\mathrm f})\\
 &=2\operatorname{Im}\!\left[
 \overline{\partial_1(w_1^{\mathrm f}+\mathrm i w_2^{\mathrm f})}\,
 \partial_\theta\partial_1(w_1^{\mathrm f}+\mathrm i w_2^{\mathrm f})\right]
 =-8\varepsilon^2g^2\|y\|^2.
\end{aligned}
\]
It follows that
\[
 \mathcal P_2\hat W=-8\varepsilon^2\|y\|^2\sin2\theta.
\]

The definition of $\hat C$ and
$\mathcal P_2\partial_{a_1}\mathcal{F}=0$ now give
\[
\begin{aligned}
 \mathcal P_2\hat C
 &=\partial_\theta^{-1}(\mathcal P_2\hat W+4\mathcal P_2\hat J)\\
 &=\partial_\theta^{-1}\!\left(-12\varepsilon^2\|y\|^2\sin2\theta\right)
 =6\varepsilon^2\|y\|^2\cos2\theta,
\end{aligned}
\]
where $\partial_\theta^{-1}\sin2\theta=-\tfrac12\cos2\theta$.
Finally, $\Delta_y^{-1}\|y\|^2=\|y\|^4/16$, so
\[
 \Delta_y^{-1}\mathcal P_2\hat J
 =-\frac{\varepsilon^2\|y\|^4}{16}\sin2\theta,\qquad
 \Delta_y^{-1}\mathcal P_2\hat C
 =\frac{3\varepsilon^2\|y\|^4}{8}\cos2\theta.
\]
Using $\langle\cos^22\theta\rangle=\langle\sin^22\theta\rangle=1/2$
and $\langle\sin2\theta\cos2\theta\rangle=0$, we find
\[
 \left\langle2\partial_{a_1}\Delta_y^{-1}\hat J\right\rangle
 =\frac{\varepsilon^2y_2\|y\|^2}{4},\qquad
 \left\langle2\partial_{a_2}\Delta_y^{-1}\hat C\right\rangle
 =\frac{3\varepsilon^2y_2\|y\|^2}{2}.
\]
Substitution in \eqref{eq:Pi3mean} yields
\[
 p_3=\frac{2\varepsilon^2}{3}
           \left(\frac14+\frac32\right)y_2\|y\|^2
     =\frac{7\varepsilon^2}{6}y_2\|y\|^2. \tag*{\qedhere}
\]
\end{proof}

\section*{Statements and Declarations}

\subsection*{Acknowledgments} 
The research of N.S. was partially supported by 
US-Japan Fusion Research Collaboration Joint Institute for Fusion Theory JPU1, and by 
JSPS KAKENHI Grant 
No. 25K07267, No. 22H04936, and No. 24K00615.  

\noindent This research was partially supported by an Australian Government Research Training Program (RTP) Scholarship at The University of Western Australia.

\subsection*{Data availability}
The Mathematica code that computes the coefficients and produces
Table~\ref{tab:visualization-residual} and
Figures~\ref{fig:mhs-pressure-surfaces}--\ref{fig:mhs-poincare}
is available at \url{https://github.com/Kooper95/An-Analytic-Toroidal-MHD-Equilibrium-with-Trivial-Euclidean-Symmetry-Group.git}.


\subsection*{AI usage statement}
The solution presented in this manuscript was discovered with the assistance of GPT-6 Astra Ultra. Generative AI tools were also used to improve the clarity and readability of the text. The authors have manually checked all results and take full responsibility for the content of the manuscript.





\bibliographystyle{alphaurl}
\bibliography{refs.bib}

\end{document}